\documentclass[10pt]{amsart} 

\usepackage{amssymb}

\usepackage{graphicx}

\newtheorem{thm}{Theorem}[section]
\newtheorem{prop}[thm]{Proposition}
\newtheorem{lemma}[thm]{Lemma}
\newtheorem{cor}[thm]{Corollary}

\theoremstyle{definition}
\newtheorem{defn}[thm]{Definition}
\newtheorem{ex}[thm]{Example}

\theoremstyle{remark}
\newtheorem{remark}[thm]{Remark}
\newtheorem{notations}[thm]{Notations}

\numberwithin{equation}{section}

\def\R{\mathbb{R}}
\def\RP{\mathbb{RP}}

\def\Z{\mathbb{Z}}

\def\P{\mathbb{P}}

\def\g{\mathfrak{g}}
\def\a{\mathfrak{a}}
\def\gl{\mathfrak{gl}}

\def\O{\mathrm{O}}

\begin{document}

\title[Deformation of surfaces in Lie sphere geometry]{Deformation and applicability of 
surfaces\\ in Lie sphere geometry}



\author{Emilio Musso}
\address{Dipartimento di Matematica Pura ed Applicata,
Universit\`a degli Studi di L'Aquila, Via Vetoio, I-67100
L'Aquila, Italy} \email{musso@univaq.it}

\author{Lorenzo Nicolodi}
\address{Di\-par\-ti\-men\-to di Ma\-te\-ma\-ti\-ca,
Uni\-ver\-si\-t\`a degli Studi di Parma, Via M. D'Azeglio 85/A,
I-43100 Parma, Italy} \email{lorenzo.nicolodi@unipr.it}

\thanks{This research was partially supported by the MIUR project
\textit{Propriet\`a Geometriche delle Variet\`a Reali e
Complesse},
and by the European Contract
Human Potential Programme, Research Training Network
HPRN-CT-2000-00101 (EDGE)}





\begin{abstract}
The theory of surfaces in Euclidean space can be naturally formulated in the more 
general context of Legendre surfaces into the space of contact elements. We address
the question of deformability of Legendre surfaces with respect to the symmetry group
of Lie sphere contact transformations from the point of view of the deformation theory
of submanifolds in homogeneous spaces. Necessary and sufficient conditions are provided for 
a Legendre surface to admit non-trivial deformations and the corresponding existence 
problem is discussed.
\end{abstract}

\maketitle

\section*{Introduction}\label{s:intro}

The classical problems of applicability of surfaces in Euclidean, projective and conformal 
geometry fit into the general theory of deformation of submanifolds in homogeneous
spaces as formulated by \'E. Cartan \cite{Ca1,Ca2} and further developed
by P. Griffiths and 
G. Jensen \cite{Gr,J}. Two submanifolds in a homogeneous space,
$f, \widetilde{f} : M \to G/K$, are $k$-th order deformations of
each other if there exists a smooth map $v : M \to G$ such that, for 
every $p\in M$, the Taylor expansions about $p$ of $\widetilde{f}$
and $v(p)\cdot{f}$ agree through $k$-th order terms; if $v$ is constant, $f$ and $\widetilde{f}$ are
congruent with respect to $G$.
Of course, for each concrete geometric situation there is a specific problem to solve. 
In Euclidean space, two surfaces are applicable in Gauss' sense
if they are first order Euclidean deformations of each other, which means they are 
isometric with respect to the induced metrics, and are congruent to second order. 
In projective $3$-space, Fubini's notion of applicability of surfaces goes to second order 
and rigidity
to third order. 
In M\"obius and Laguerre geometry, second
order deformable surfaces coincide with isothermic and $L$-isothermic surfaces, 
respectively \cite{M1, MN2}.

\smallskip
This paper studies the deformation problem for surfaces in another classical
geometry: Lie sphere geometry. It is the outcome of our attempts to understand 
Lie-applicable surfaces within the general theory of deformation.
Lie-applicable surfaces were considered by Blaschke and his
collaborators in the 1920s \cite{Bl2}. 
Recently, interest in Lie-applicable surfaces has reappeared in the work of 
Ferapontov concerning the relation between Lie sphere geometry of
hypersurfaces and the theory of integrable systems \cite{Fe1,Fe3}.

\smallskip
To put our discussion in perspective let us recall some facts about Lie sphere geometry. 
Any smooth immersion of an oriented surface into $3$-space has a contact lift to the unit sphere bundle
$\Lambda$ of $S^3=\R^3\cup\{\infty\}$ called the Legendre lift. The unit sphere bundle is
acted on transitively by the group of Lie sphere transformations. 
This group is the group of contact transformations generated by the conformal 
transformations of $S^3$ together with the group of
normal shifts, which transform
an oriented surface to its parallel surface at a fixed oriented distance in
the normal direction. Thus the Lie sphere group acts on the set of Legendre surfaces.
The principal aim of Lie sphere geometry is to study the properties of an immersion which are
invariant under this action \cite{Bl2, Ce, Pi}.
 Blaschke associated with any Legendre surface  
a canonical coframe $\mathcal{C}$ and
proved that the position of a generic Legendre surface is 
completely determined by its canonical coframe, up to Lie sphere transformations.
However, he observed that there are exceptions to this rigidity result. 
Accordingly, two non-congruent
Legendre immersions $f,\widetilde{f}$ are
called \textit{Lie-applicable} if
$\mathcal{C} = \widetilde{\mathcal{C}}$. 
As we will see in Section \ref{s: defo}, examples of Lie-applicable surfaces include the Legendrian 
lifts of isothermic and $L$-isothermic surfaces, which are known to constitute
integrable systems \cite{Bu, Be, MN4}.

\smallskip

In the paper, we think of $\Lambda$ as a homogeneous space 
of the identity component $\widetilde{G}$ 
of the Lie sphere group. We will prove that two (nondegenerate) Legendre immersions 
are Lie-applicable if and only if they are second order deformations of
each other; that two Legendre immersions are always local first order deformations of each 
other; and that they have third order rigidity. Further, we show how to recognize that a Legendre
surface is deformable and in this case how to find all its deformations.

\smallskip
In Section \ref{s: basic}, we collect some background material about the Lie sphere
geometry of surfaces and develop the
method of moving frames in this context  (see \cite{Bl2, Ce}). 
We identify $\Lambda$ with the space of isotropic 2-spaces 
in $\R^{4,2}$ and linearize the action of Lie transformations. In this model,
the identity component of $\O(4,2)$ acts on $\Lambda$ by contact
diffeomorphisms and provides a double cover of $\widetilde{G}$.   
We then apply the method of moving frames to study Legendre surfaces and
recall the construction of a canonical lift to the group $\widetilde{G}$ 
under a natural nondegeneracy assumption.
For any Legendre surface, we introduce the canonical coframe $(\alpha^1,\alpha^2)$, 
which turns out to be the analogue of that considered by Blaschke \cite{Bl2}, 
and the quadratic and cubic forms of the surface. 
We define a set of local differential invariants $q_1,q_2,p_1,p_2,r_1,r_2$ for
a Legendre surface and
relate them to the classical ones discussed by Blaschke and Ferapontov.
We then deduce the compatibility conditions, 
which play the role of the Gauss--Codazzi equations for a Legendre surface in Lie 
sphere geometry. 
The functions $q_1,q_2,p_1,p_2$ are completely determined by the canonical coframe, 
while $r_1,r_2$ govern the extrinsic geometry of the Legendre immersion.

\smallskip
In Section \ref{s: defo}, we investigate the class of
Legendre surfaces which are not determined by the canonical coframe. We 
take the point of view of the deformation theory of submanifolds in homogeneous spaces.
We introduce the concept of deformation and discuss the related questions of analytic 
contact and applicability. We study the problem of second order deformation
of Legendre immersions and prove that two nondegenerate Legendre immersions are
second order deformations of each other precisely when they are Lie-applicable, or equivalently,
when they have the same 
quotient of cubic to quadratic forms (see Theorem \ref{main-theorem}).

\smallskip
In Section \ref{s: infinitesimal-defo}, for a given Legendre immersion,
we introduce a suitable linear connection 
on the trivial bundle $M^2 \times \R^3$ and
show that the existence of 
non-zero parallel sections with respect to this connection is a necessary and 
sufficient condition for the Legendre immersion to have non-trivial deformations. 
In particular, we show that the non-trivial
deformations of nondegenerate Legendre immersion depend on three parameters at most.
The above characterization allows the
introduction of free parameters in the Maurer-Cartan form of the canonical
frame of a deformable surface without violating the structure
equations, and suggests the existence of a B\"acklund
transformation for the class of deformable Legendre surfaces. We will return on
this topic elsewhere.


\smallskip
In Section \ref{s: examples}, we discuss some examples of deformable Legendre immersions.

\smallskip
In the last section, we study the existence of
deformations. 
We use the characterization of deformable Legendre immersions
in terms of parallel sections
to 
set up the exterior differential system of a deformation. We then prove 
that this system is in involution in Cartan's sense and that its general 
solutions depend on six functions in one variable.

\bigskip

\section{Surface theory in Lie sphere geometry}\label{s: basic}

In this section, we briefly recall the basic structure of Lie sphere geometry
and develop the method of moving frames for immersed surfaces in the context of that
geometry. More details about Lie sphere geometry are given in the recent monograph 
of Cecil \cite{Ce}, in the book of Blaschke \cite{Bl2}, or in Lie's original work
\cite{LS}.

\subsection{ Legendre immersions}

Let $S^3$ be the unit sphere in $\R^4$ and identify the unit tangent bundle $\Lambda = T_1S^3$
with the set of all pairs $(v,\xi)\in S^3\times S^3$ such that
$v$ is orthogonal to $\xi$, i.e.,
\[
  \Lambda=\{(v,\xi)\in S^3\times S^3 \subset \R^4\times \R^4 \, |  \, v\cdot
   \xi=0\}.
    \]
Let $\pi_1, \, \pi_2 : T_1S^3 \to S^3$ denote the restrictions to $\Lambda$ of the canonical 
projections of $S^3 \times S^3$ onto its factors. Then the equation
$d\pi_1 \cdot \pi_2 =0$ defines a $4$-dimensional contact distribution $\mathcal D$
on $\Lambda$. 

\medskip
If $F : M^2 \to S^3$ is an immersed surface oriented by a field of 
unit normals $n$, then $(F,n) : M^2 \to \Lambda$ is an integral submanifold of $\mathcal D$.
In general, an immersion $f : M^2 \to \Lambda$ which is an integral submanifold of $\mathcal D$
is called a \textit{Legendre surface}. 
The \textit{Lie sphere group}, that is, the group generated by the conformal transformations 
of $S^3$ together with the group of normal shifts, which transform
an oriented surface to its parallel surface at a fixed oriented distance in
the normal direction, preserves the contact distribution $\mathcal D$ and acts naturally on 
the space of Legendre surfaces. 

\medskip
If $f= (F,n) : M^2 \to \Lambda$ is a Legendre immersion, the smooth map $F : M^2 \to S^3$
need not be an immersion. However, without loss of generality, 
we can always assume (applying if necessary a normal shift) that $F$ is locally an immersion. 
This follows from a result of Pinkall 
\cite{Pi} (see also \cite{Ce}) stating that, for each $p\in M^2$, there exists $t\in [0,\pi)$ 
for which the parallel
surface $F_t = \cos t \,F + \sin t \, n$ is locally an immersion.

We recall that the \textit{curvature sphere} at $p\in M^2$ corresponding to a principal 
curvature $k_i$ is the oriented sphere in oriented contact with $F(M^2)$ at $F(p)$ and 
centered at the focal point determined by the principal curvature $k_i$. The notion of a
curvature sphere is invariant under Lie transformations.
For instance, if $F : M^2 \to \R^3$ is an immersed surface oriented by the field of normals 
$n : M^2 \to S^2$, then $(F,n) : M^2 \to \R^3\times S^2 \subset \Lambda$ is a Legendre
surface. In this case, the curvature spheres at a point $p\in M^2$ are the oriented spheres 
$\sigma_i(p)$ centered at
$F(p) + k_i(p)^{-1}{n}(p)$, with signed radius
$k_i(p)^{-1}$, $i=1,2$. When $k_i(p)=0$, then $\sigma_i$ is the
oriented tangent plane of the surface at $F(p)$. 

\begin{defn}[Nondegeneracy condition] 
We say that a Legendre surface $f= (F,n) : M^2 \to \Lambda$ is \textit{nondegenerate}  
if $F$ is umbilic free and both the {curvature spheres} corresponding to the principal
curvatures $k_1,\, k_2$ are immersions into the space of oriented $2$-spheres in $S^3$ 
(including points). 
\end{defn}

\subsection{Moving Lie frames for Legendre surfaces}

Let $\R^{4,2}$ denote $\R^6$ with the symmetric bilinear form
\begin{equation}\label{scalar-product}
 \langle X,Y\rangle =
   -(x^0y^5+x^5y^0)-(x^1y^4+x^4y^1)+x^2y^2+x^3y^3={^tX} gY
    \end{equation}
of signature $(4,2)$,
where $(x^J)$ and $(y^J)$ are the coordinates of $X$ and $Y$ with
respect to the standard basis $(\epsilon_0,\dots,\epsilon_5)$ of $\R^6$.
Let $G$ be the connected component of the identity of the group 
\[
 \{A \in \mathrm{GL}(6,\R) \, |\, {^tA}gA=g \}
  \]
of linear transformations which leave the form (\ref{scalar-product}) invariant
and let
$\mathfrak{g} = \{B\in \mathfrak{gl}(6,\R) \, |\, {^tB}g+gB=0 \}$ be
its Lie algebra. For each
$A\in G$, we denote by $A_J=A\cdot \epsilon_J$ the $J$-th column
vector of $A$. Regarding each of the vectors $A_J$ as a vector-valued function $v:G\to \R^{6}$
on $G$, since the $A_J$ form a basis, there exist unique 1-forms 
$\omega^I_J$, $I,J\in\{0,1,\dots,5\}$, so that
\begin{equation}\label{dA}
  dA_J=\omega^I_JA_I,\quad J=0,\dots,5.
   \end{equation}
(We use the summation convention on repeated indices.)
The $1$-forms $\omega^I_J$ are the components of the left-invariant
{Maurer--Cartan form} $\omega=A^{-1}dA$ of $G$. They are connected by relations
obtained from the differentiation of $\langle A_I,A_J\rangle = 
g_{IJ}$, $I,J\in\{0,1,\dots,5\}$, which are
\begin{equation}\label{omega-g}
 {^t\omega} g + g\omega = 0, \quad \text{or} \quad \omega^K_Ig_{KJ} + \omega^K_Jg_{KI}, 
  \quad I,J\in\{0,1,\dots,5\},
   \end{equation} 
and reflect the structure of the Lie algebra $\mathfrak{g}$.
The forms
\[
 \omega^0_0,\,\omega^1_1,\,\omega^0_1,\,\omega^1_0,\,\omega^2_0,\,\omega^3_0,\,\omega^2_1,\,
  \omega^3_1,\,\omega^4_0,\,\omega^3_2,\,\omega^0_2,\,
   \omega^0_3,\,\omega^1_2,\,\omega^1_3,\,\omega^0_4
    \]
yield a left-invariant coframe field on $G$.
Differentiating (\ref{dA}), we obtain the structure equations of $G$, 
which are
\begin{equation}\label{d-omega}
 d\omega = -\omega\wedge \omega, \quad \text{or} \quad d\omega^I_J = -\omega^I_K\wedge\omega^K_J,
  \quad I,J\in\{0,1,\dots,5\}.
   \end{equation} 
For each  $X\in G$, the Maurer--Cartan form $\omega$ transforms as follows
\begin{equation}\label{1.1.6}
  R_X^\ast(\omega) = X^{-1}\omega X.
   \end{equation}

\medskip
The projectivization $\mathcal{Q} = \P(\mathcal{L})$ of 
the light cone $\mathcal{L}$ of $\R^{4,2}$ is known in the classical literature as the \textit{Lie quadric}. 
In Lie sphere geometry, the Lie quadric parametrizes the set of all oriented 2-sphere in $S^3$, 
including points, and the
lines in $\mathcal{Q}$ correspond to parabolic pencils of spheres in oriented contact.
The set of all lines in $\mathcal{Q}$, that is, the isotropic Grassmannian of null 2-planes through 
the origin in $\R^{4,2}$, forms a smooth manifold which can be identified with $\Lambda$
(for more details see \cite{Ce}).
Under this identification, the group $G$ acts transitively on
$\Lambda$ by the usual action of $G$ on the Grassmannian $G_2(\R^{4,2})$ and
preserves the contact structure.
The projection map
\begin{equation}\label{fibration}
 \pi_{\Lambda} : A\in G\to A[\epsilon_0\wedge\epsilon_1]=[A_0\wedge A_1]\in
  \Lambda = G/G_0,
   \end{equation}
defines a principal $G_0$-bundle over $\Lambda$, where $G_0$ is the isotropy subgroup at the 
chosen origin $[\epsilon_0\wedge\epsilon_1]$. 

\medskip
A Legendre surface $f : M^2 \to \Lambda$ may then be represented by two maps $F_0 , F_1 : M^2\to
\mathcal{L}$ such that $f=[F_0\wedge F_1]$, $\langle dF_0,F_1\rangle = 0$ and $\langle F_0,F_1\rangle = 0$.
Of course, such a representation is not unique.
For example, if $F : M^2 \to \R^3$ is any smooth immersion in $\R^3$, 
oriented by a field of unit normals
$n : M^2\to S^2$, then the \textit{Legendre lift} $f=[F_0\wedge F_1]$ is given by
\begin{equation}\label{2.1}
 \left\{\begin{array}{ll}
  F_0={^t\left(1,\frac{1}{\sqrt{2}}F^1,F^2,F^3,-\frac{1}{\sqrt{2}}F^1,\frac{1}{2}F\cdot
    F \right)},\\
     F_1=\frac{1}{\sqrt{2}}{^t\left(0,\frac{1}{2}(1+n^1),n^2,n^3,
      \frac{1}{\sqrt{2}}(1-n^1),n\cdot F\right)}.
       \end{array}\right. 
        \end{equation}

\medskip

A \textit{frame field} along a Legendre surface $f:M^2 \to \Lambda$ is a smooth map
$A : U \to G$ defined on some open subset of $M^2$ such that $f=[A_0\wedge A_1]$.
For each local frame $A : U\to G$ we let
$\alpha= A^\ast \omega$. 
The Legendre condition simply means that the form $\alpha^4_0$ vanishes
identically.
Any other local frame is given by $\widetilde{A}=A\cdot X$, for some smooth map $X : U\to G_0$,
and the $1$-forms $\alpha$ and $\widetilde{\alpha}$ are related by
\begin{equation}\label{2.1.1}
 \widetilde{\alpha}=X^{-1}dX+X^{-1}\alpha X.
  \end{equation} 
The totality of frames along $f$ is the principal $G_0$-bundle $\mathcal{F}_0(f) \to M^2$, where
\[
  \mathcal{F}_0(f)=\{(p,A)\in
    M\times G \mid f(p)=[A_0\wedge A_1]\}.
     \]

\medskip
\textit{The canonical frame.} 
Following the usual practice in the method of moving frames, we can construct a canonical lift to
the group $G/\Z_2$ for any nondegenerate Legendre surface. The idea of the procedure is at each step
to normalize the Maurer--Cartan matrix of a frame along $f$ as much as possible, and then 
take the exterior derivative of the equations expressing this normalization, thereby leading to 
the next step. Similar preferred frames have been given by Blaschke 
\cite{Bl2} and Ferapontov \cite{Fe1, Fe3}. Here, we skip the construction.

\begin{thm}[Existence of the canonical frame]
Let $f : M^2\to \Lambda$ be a nondegenerate Legendre immersion of an oriented surface $M^2$. Then 
there exists a unique 
lift $[A]:M^2\to G/\Z_2$ satisfying the Pfaffian equations
\begin{equation}\label{2.3.2}
 \alpha^4_0=\alpha^2_0=\alpha^3_1=
  \alpha^3_2=\alpha^1_0-\alpha^2_1=\alpha^0_1-\alpha^3_0=\alpha^0_2=\alpha^1_3=0
   \end{equation}
with the independence condition
\begin{equation}\label{2.3.3}
  \alpha^3_0\wedge \alpha^2_1>0. 
   \end{equation}
$[A]$ is called the {\em canonical frame field} along $f$. Let $\alpha^1 = \alpha^3_0$ and $\alpha^2 =\alpha^2_1$. 
The coframe $(\alpha^1, \alpha^2)$ on $M^2$ is referred to as the {\em canonical coframe} of $f$.
The bundle of canonical frames along $f$ we will be denoted by $\mathcal{F}(f) \to M^2$ .
\end{thm}

\begin{defn}
The {\em quadratic form} $\Phi$ and the {\em cubic form} $\Psi$ of the immersion $f$ are defined by
\begin{equation}\label{forms}
 \Phi = - \alpha^1\alpha^2, \qquad \Psi = -(\alpha^1)^3 + (\alpha^2)^3,
  \end{equation}
respectively.
The quotient $\mathcal{B}=\Psi/\Phi$ of the cubic form $\Psi$ to the quadratic form $\Phi$ is a 
well defined map $\mathcal{B}:T(M)\to \R \cup \{\infty\}$ which is a rational function when restricted to
the tangent spaces $T_p(M)$. $\mathcal{B}$ is called the {\em Fubini--Blaschke invariant} of $f$.
\end{defn}

\medskip
\textit{The invariant functions.} 
The only non zero components of the Maurer--Cartan form $\alpha$ of the canonical frame are
$\alpha^1 = \alpha^3_0$, $\alpha^2= \alpha^2_1$, and $\alpha^0_0$, $\alpha^1_1$, $\alpha^1_2$,
$\alpha^1_1$, $\alpha^0_3$, $\alpha^0_4$.
From the exterior differentiation of (\ref{2.3.2}) and the
structure equations, it follows that there exist smooth
functions $q_1,q_2$, $p_1,p_2$, and $r_1,r_2$ such that
\begin{equation}\label{2.3.4}
 \left\{ \begin{array}{lll} 
 \alpha^0_0=-2q_1\alpha^1 + q_2\alpha^2, &
  \alpha^1_1= -q_1\alpha^1+2q_2\alpha^2,\\
   \alpha^0_3= r_1\alpha^1+p_2\alpha^2,&
    \alpha^1_2= p_1\alpha^1 + r_2 \alpha^2,\\
     \alpha^0_4=-r_2\alpha^1 + r_1\alpha^2,
      \end{array}\right. 
       \end{equation}
We shall refer to these functions as the \textit{invariant functions} of $f$. 
Using the structure equations,
we obtain
\begin{equation}\label{2.3.5}
  d\alpha^1=\alpha^0_0\wedge\alpha^1,\quad
   d\alpha^2=\alpha^1_1\wedge\alpha^2,
    \end{equation}
\begin{equation} \label{2.3.6}
   d\alpha^0_0=(\alpha^2-\alpha^0_3)\wedge\alpha^1,\quad
    d\alpha^1_1=(\alpha^1-\alpha^1_2)\wedge\alpha^2,
     \end{equation}
\begin{equation}\label{2.3.7}
 d\alpha^1_2=-\alpha^1_1\wedge\alpha^1_2,\quad
  d\alpha^0_3=-\alpha^0_0\wedge\alpha^0_3,\quad
   d\alpha^0_4=-(\alpha^0_0+\alpha^1_1)\wedge\alpha^0_4.
    \end{equation}
In terms of the invariant functions, these equations become
  \begin{equation}\label{2.3.8}
   d\alpha^1=-q_2\alpha^1\wedge \alpha^2,\quad
    d\alpha^2=-q_1\alpha^1\wedge \alpha^2 
     \end{equation}
\begin{equation}\label{2.3.9}
  \left\{ \begin{array}{rcl} 
 -2dq_1\wedge \alpha^1+dq_2\wedge \alpha^2 &=& (p_2-q_1q_2-1)\alpha^1\wedge\alpha^2,\\ 
   -dq_1\wedge \alpha^1+2dq_2\wedge \alpha^2 &=&(-p_1+q_1q_2+1)\alpha^1\wedge \alpha^2,
     \end{array}\right. 
      \end{equation}
\begin{equation}\label{2.3.10}
  \left\{ \begin{array}{rcl} 
   dr_1\wedge \alpha^1 + dp_2\wedge \alpha^2 &=& (2q_2r_1+3q_1p_2)\alpha^1\wedge \alpha^2,\\ 
    dp_1\wedge \alpha^1+ dr_2\wedge \alpha^2 &=& (2q_1r_2+3q_2p_1)\alpha^1\wedge \alpha^2,\\ 
     -dr_2\wedge \alpha^1+ dr_1\wedge \alpha^2 &=& 4(q_1r_1-q_2r_2)\alpha^1\wedge \alpha^2.
       \end{array}\right. 
        \end{equation}
Equations (\ref{2.3.8}) and (\ref{2.3.9}) tell us that the invariant functions 
$q_1$, $q_2$, $p_1$, and $p_2$ are determined by the canonical coframe. 
Equations (\ref{2.3.10})
can be viewed as compatibility conditions arising from the
fact that the canonical coframe is obtained from the Legendre
immersion. Thus, we may think of (\ref{2.3.8}) and (\ref{2.3.9})
as the Gauss equations and of (\ref{2.3.10}) as 
the Codazzi-Mainardi equations of the immersion.

\medskip
\textit{Relations with Euclidean geometry.} 
Let $f : M^2 \to \Lambda$ be the Legendre lift of an oriented immersion
$F : M^2 \to \R^3$ with Gauss map ${n}= (n^1,n^2,n^3)$ and suppose $f$ is nondegenerate.
Let $(u,v)$ be curvature line coordinates on $M^2$.
Then the 
canonical coframe takes the form
\[  
 \alpha^1=\frac{1}{k_1-k_2}\left(\sqrt{{e}^{-1}{g}}\, (k_1)_u\,(k_2)_v^2\right)^{\frac{1}{3}}dv,\quad
   \alpha^2=\frac{1}{k_2-k_1}\left(\sqrt{{e}{g}^{-1}}\, (k_1)_u^2\,(k_2)_v\right)^{\frac{1}{3}}du,
   \] 
where $e$ and $g$ are the coefficients of the first fundamental
form of $F$ with respect to the coordinate system $(u,v)$.  Moreover,
the quadratic and cubic forms takes the form
\[
   \Phi= \frac{1}{(k_1-k_2)^2}(k_1)_u(k_2)_v\, {du}\,{dv}
      \]
\[
 \Psi= - \frac{(k_1)_u\, (k_2)_v }{(k_1-k_2)^3 \sqrt{{e}{g}}} \left({e}\, (k_1)_u  \,du^3
     +{g}\, (k_2)_v \,dv^3\right).  
    \] 
Observe that $\alpha^1=(\beta\gamma^2)^{1/3}dv$ and
$\alpha^2=(\beta^2\gamma)^{1/3}du$, where $\beta$ and $\gamma$
are given by
\[
 \beta =\frac{1}{k_1-k_2} \sqrt{eg^{-1}}\, (k_1)_u, \quad
   \gamma =\frac{1}{k_2-k_1} \sqrt{e^{-1}g}\, (k_2)_v.
    \] 
Using the above structure equations, the invariants $q_1, q_2, p_1, p_2$ can be expressed in 
terms of $\beta$ and $\gamma$. For example,
\[
 q_1= -\frac{1}{3(\beta^2\gamma)^{2/3}}\left(2\beta_v+\frac{\beta}{\gamma}\gamma_v\right),
 \quad
  q_2= \frac{1}{3(\beta\gamma^2)^{2/3}}\left(\frac{\gamma}{\beta}\beta_u+2\gamma_u\right).
   \]
The invariants $\beta$ and $\gamma$ correspond to the invariants $p$ and $q$ considered by 
Ferapontov \cite{Fe1}, p. 207, and Blaschke \cite{Bl2} in the construction of the canonical frame.
Note that the vanishing of both $\beta$ and $\gamma$ is equivalent to the
condition that the principal curvatures are constant
along the corresponding principal foliations. This property characterizes the \textit{cyclides of Dupin}.
If one of the two principal curvatures is constant along the corresponding
principal foliation, the surface is the envelope of a one-parameter 
family of oriented spheres (including planes and
point-spheres), and we are in presence of a canal surface. 

\medskip
\begin{remark}\label{rk: 3-web}
Associated with any nondegenerate Legendre immersion $f$ there is the $3$-web formed 
by the \textit{asymptotic lines}
of the quadratic form $\Phi$ and by the \textit{cyclidic curves}\footnote{The family of
curves which are orthogonal to the cyclidic curves with respect to
the quadratic form $\Phi$ is called \textit{anti-cyclidic system}.}, i.e., 
the curves along the zero-directions of the cubic form
$\Psi$. In view of (\ref{forms}), the
curves of this web can be defined in terms of the canonical
coframe by the Pfaffian equations
\begin{equation}\label{2.3.1.1}
 \alpha^1=0,\quad \alpha^2=0,\quad
  \alpha^1-\alpha^2=0.
   \end{equation} 
The connection form of the $3$-web is the $1$-form $\zeta_w$,
uniquely determined by the equations
\[
  d\alpha^1=\zeta^w\wedge\alpha^1,\quad d\alpha^2=\zeta^w\wedge
  \alpha^2.
   \] 
Equations (\ref{2.3.4}) and (\ref{2.3.5}) yield
\begin{equation}\label{2.3.9.1}
 \zeta_w = -q_1\alpha^1+q_2\alpha^2.
  \end{equation}
From (\ref{2.3.9.1}) one then computes the curvature of the $3$-web, which is
\begin{equation}\label{2.3.9.2}
  R_w=\frac{1}{3}(p_2-p_1).
   \end{equation}
The surfaces for which the curvature vanishes identically are called diagonally
cyclidic (see \cite{Bl2, Fe1}).

\end{remark}

\bigskip

\section{Contact, deformation and applicability}\label{s: defo}

Let us recall the general notion of deformation \cite{Gr,J}.

\begin{defn}
Let $G/K$ be a homogeneous space and let $f,\tilde{f} : N\to
G/K$ be two smooth maps. We say that $f$ and $\tilde{f}$ are
\textit{$k$-th order deformations of each other} if there exists
a smooth map ${B}:N\to G$ such that, for each point
$p\in N$, $\tilde{f}$ and
${B}(p)f$ have the same $k$-th order jets at $p$, i.e., they have
analytic contact of second order at $p$. The map ${B}$ is 
said to be a \textit{$k$-th order
deformation}. When ${B}$ is constant the deformation
is said to be \textit{trivial}. 
A map $f:N\to G/K$ is said to be \textit{deformable of order $k$} 
if it admits a non-trivial $k$-th order deformation.
\end{defn}

First, we will express the condition of analytic contact in the special case 
of mappings from a 2-dimensional manifold $M$ into $\Lambda$.  For this we 
need to introduce some notations.

\subsection{Analytic Contact}\label{ss: 3.1}

Let $(x^1,x^2)$ be a local coordinate system on an open set $U$
of $M$. Let $E$ be a vector space and let $S^h(U)\otimes E$ 
denote the symmetric $E$-valued $k$-forms on $U$. The symmetric tensor
product of $s\in S^h(U)$ and $t\in S^k(U)$ will be denoted by $s\cdot t$. 
An element $T$ of $S^h(U)\otimes E$ has a local expression 
$$
T=T_{i_1\dots i_h}dx^{i_1}\dots dx^{i_h},
$$
where the coefficients $T_{i_1\dots i_h}$ are $E$-valued
smooth maps, which are totally symmetric in the indices $i_1,\dots,i_h$. We
then define the $k$-th order derivative of $T$ to be the
$E$-valued symmetric form of order $h+k$ given by
$$
\delta^k(T)=\frac{\partial^k T_{i_1 \dots i_h}}{\partial
x^{i_{h+1}}\dots\partial
x^{i_{h+k}}}dx^{i_1}\dots dx^{i_h}dx^{i_{h+1}}\dots dx^{i_{h+k}}.
$$
The definition depends on the choice of the local
coordinates. Given a pair $T_0,T_1\in S^h(U)\otimes E$ of
$E$-valued symmetric $h$-forms and a $2\times 2$ matrix
$\rho=(\rho^a_b)\in S^k(U)\otimes \gl(2,\R)$ of symmetric
$k$-forms we set
$$
 (T_0,T_1)\rho = (\rho^0_0T_1+\rho^1_0 T_1,\rho^0_1
  T_0+\rho^1_1 T_1).
$$
We can state the following

\begin{lemma}\label{lemma1}
Let $f=[F_0\wedge F_1]: M^2 \to \Lambda$ and
$\tilde{f}=[\tilde{F}_0\wedge \tilde{F}_1]:M\to \Lambda$ be two
smooth maps. Then, $f$ and $\tilde{f}$ agree to second order at $p\in M$, i.e.,
they have the same second
order jets at $p$, if and only if, for every local
coordinate system $(x^1,x^2)$ about $p$, there
exist
$$
 \rho_0\in \gl(2,\R),\quad \rho_1\in T^*(M)_p\otimes
  \gl(2,\R),\quad \rho_2\in S^2(M)|_p\otimes \gl(2,\R)
$$
\noindent such that
 \begin{equation}\label{3.1}
   \begin{aligned}
    &(\tilde{F}_0|_p,\tilde{F}_1|_p)=(F_0|_p,F_1|_p)\rho_0,\\ 
     &(\delta\tilde{F}_0|_p,\delta \tilde{F}_1|_p)=(\delta F_0|_p,\delta
      F_1|_p)\rho_0+(F_0|_p,F_1|_p)\rho_1,\\ 
       &(\delta^2\tilde{F}_0|_p,\delta^2 \tilde{F}_1|_p)=(\delta^2 F_0|_p,\delta^2
        F_1|_p)\rho_0+(\delta F_0|_p,\delta
         F_1|_p)\rho_1+(F_0|_p,F_1|_p)\rho_2,
         \end{aligned}
          \end{equation}
\end{lemma}

\begin{proof} 
Let $(x^1,x^2)$ be a local coordinate system  on
an open neighborhood $U$ of $p$. As $G$ acts
transitively on $\Lambda$, we may assume that
\[
f(p)=\tilde{f}(p)=[\epsilon_0\wedge \epsilon_1].
\]
The map
\[
  y=(y^1,\dots,y^5)\in \R^5 \mapsto [X_0(y)\wedge X_1(y)] \in \Lambda
  \]
defined by
\begin{equation}
 \left\{\begin{array}{l}
   X_0(y) = {^t\left(1,0,y^1,y^2,y^3,\frac{1}{2}[(y^1)^2+(y^2)^2]\right)},\\
    X_1(y) = {^t\left(0,1,y^4,y^5,\frac{1}{2}[(y^4)^2+(y^5)^2],y^1y^4+y^2y^5-y^3\right)},
     \end{array}\right.
      \end{equation}
is a local coordinate system of $\Lambda$ centered at
$[\epsilon_0\wedge \epsilon_1]$. Then, there exists an open
neighborhood  $U'\subset U$ of $p$ and smooth maps
$h,\tilde{h} : U'\to \R^5$ such that
\[
 f|_{U'}=[(X_0\circ h)\wedge (X_1\circ h)],\quad
  \tilde{f}|_{U'}=[(X_0\circ \tilde{h})\wedge (X_1\circ
   \tilde{h})].
    \]
Thus, $f$ and $\tilde{f}$ have second order analytic
contact at $p$ if and only if the maps $G_a=X_a\circ h$ and $
\tilde{G}_a=X_a\circ \tilde{h}$, $a=0,1$, satisfy
\begin{equation}\label{3.2}
  G_a(p)=\tilde{G}_a(p),\quad \delta G_a|_p=\delta
   \tilde{G}_a|_p,\quad \delta^2 G_a|_p=\delta^2 \tilde{G}_a|_p,\quad
     a=1,2.
     \end{equation}
 Let us write
\begin{equation}\label{3.3}
 (F_0,F_1)=(G_1,G_2)a,\quad
  (\tilde{F}_0,\tilde{F}_2)=(\tilde{G}_0,\tilde{G}_2)b.
\end{equation}
\noindent where $a,b:U'\to GL(2,\R)$ are smooth maps. Using
(\ref{3.2}) and (\ref{3.3}), a direct computation shows that this
is equivalent to (\ref{3.1}), where $\rho_0$, $\rho_1$ and
$\rho_2$ are given by
\begin{equation}\label{3.4}
 \left\{ \begin{array}{lll} \rho_0=a(p)b(p)^{-1},\\ \rho_1=(\delta
  a|_p-\rho_0 \delta b|_p)b(p)^{-1},\\
   \rho_2=\left(\delta^2a|_p-\rho_0\delta^2b|_p-2\rho_1 \delta
    b|_p\right)b(p)^{-1}.
     \end{array}\right. 
\end{equation}
\end{proof}

\begin{remark}
From the proof of this lemma we see that
$f$ and $\tilde{f}$ have first order analytic contact at $p$ if
and only if there exist $\rho_0\in \gl(2,\R)$ and $\rho_1\in
T^*(M)_p\otimes \gl(2,\R)$ such that
\[
 \left\{ \begin{array}{l}
  (\tilde{F}_0|_p,\tilde{F}_1|_p)=(F_0|_p,F_1|_p)\rho_0,\\ (\delta
   \tilde{F}_0|_p,\delta \tilde{F}_1|_p)=(\delta F_0|_p,\delta
     F_1|_p)\rho_0+(F_0|_p,F_1|_p)\rho_1.
     \end{array}\right.
      \]
\end{remark}

\subsection{Deformation of Legendre surfaces}\label{ss: 3.2}

\begin{notations}
Given two maps $f,\tilde{f} : M^2 \to \Lambda$, let $\mathcal{F}_0(f)$ and
by $\mathcal{F}_0(\tilde{f})$ be the $G_0$-bundles induced on $M$ by
$f$ and $\tilde{f}$, respectively. We let $j:\mathcal{F}_0(f) \to G$
and $\tilde{j}:\mathcal{F}_0(\tilde{f})\to G$ be the natural maps
\[
  j:(p,A)\in \mathcal{F}_0(f)\to A\in G,\quad
   \tilde{j}:(p,\tilde{A})\in \mathcal{F}_0(\tilde{f})\to \tilde{A}\in G.
    \]
The pull back of the Maurer-Cartan form of $G$ by $j$
and $\tilde{j}$ will be denoted by $\omega=(\omega^I_J)$ and by
$\tilde{\omega}=(\tilde{\omega}^I_J)$, respectively. If $A:U\to G$
and $\tilde{A}:U\to G$ are local cross sections of
$\mathcal{F}_0(f)$ and $\mathcal{F}_0(\tilde{f})$, respectively, then
the $\g$-valued $1$-forms $A^{-1}dA$ and
$\tilde{A}^{-1}d\tilde{A}$ will be denoted by $\alpha$ and
$\tilde{\alpha}$.
\end{notations}

\subsubsection{Deformations of order zero} 

A deformation of order zero between $f$ and $\tilde{f}$ is a smooth
map $B:M\to G$ such that $\tilde{f}(p)=B(p)f(p)$, for every $p\in
M$. Thus, $B$ induces a bundle isomorphism $\mathcal{B}:
\mathcal{F}_0(f)\to \mathcal{F}_0(\tilde{f})$ defined by the formula
\[
 \mathcal{B}:(p,A)\in \mathcal{F}_0(f)\to
  (p,B(p)A)\in \mathcal{F}_0(\tilde{f}),\quad \forall \, (p,A)\in
   \mathcal{F}_0(f).
    \]
Conversely, every bundle isomorphism between
$\mathcal{F}_0(f)$ and $\mathcal{F}_0(\tilde{f})$ arises from a deformation
of order zero.
\bigskip

\subsubsection{First order deformations}

\begin{thm}\label{teorema2}
A zero-th order deformation $B:M\to G$ of the two
maps $f,\tilde{f}:M\to \Lambda$ is of first order if and only
if the bundle isomorphism $\mathcal{B} : \mathcal{F}_0(f)\to
\mathcal{F}_0(\tilde{f})$ satisfies
\begin{equation}
 \tilde{\omega}^2_0=\mathcal{B}^*(\omega^2_0),\quad  
  \tilde{\omega}^3_0=\mathcal{B}^*(\omega^3_0),\quad
   \tilde{\omega}^4_0=\mathcal{B}^*(\omega^4_0)\quad
    \tilde{\omega}^2_1=\mathcal{B}^*(\omega^2_1) \quad
     \tilde{\omega}^3_1=\mathcal{B}^*(\omega^3_1). \label{3.5}
      \end{equation}
\end{thm}

\begin{proof} 
We have to show that for every local cross section $A:U\to G$ of
 $\mathcal{F}_0(f)$, the cross section $\tilde{A} = BA : p\in U\to B(p)A(p)\in G$
of $\mathcal{F}_0(\tilde{f})$ satisfies
\[
 \tilde{\alpha}^2_0=\alpha^2_0,\quad  \tilde{\alpha}^3_0=\alpha^3_0,\quad
  \tilde{\alpha}^4_0=\alpha^4_0\quad
   \tilde{\alpha}^2_1=\alpha^2_1 \quad
    \tilde{\alpha}^3_1=\alpha^3_1.
     \]
Recall that $B$ is a first order deformation if and only if the two maps 
$\tilde{f}$ and $B(p)f$ have first order analytic contact at $p$, for
each point $p\in M$. 
The map $A' = B(p)A : U \to G$ is a frame
along $B(p)f$ and
$B$ is a first order deformation if and only if the maps
\[
  F':q\in U\to [A'_0|_q\wedge A'_1|_q]\in \Lambda, \quad
  F:q\in U\to [\tilde{A}_0|_q\wedge \tilde{A}_1|_q]\in \Lambda
   \]
have first order analytic contact at $p$. 
From the
characterization of analytic contact, it follows that $F$ and $F'$
have first order analytic contact at $p$ if and only if there
exist $\rho_0(p)\in \gl(2,\R)$ and $\rho_1(p)\in T^*(M)_p\otimes
\gl(2,\R)$ such that
\begin{equation}\label{3.6}
 \left\{ \begin{array}{ll}
  (\tilde{A}_0|_p,\tilde{A}_1|_p)=(A'_0|_p,A'_1|_p)\rho_0(p),\\
   (\delta \tilde{A}_0|_p, \delta \tilde{A}_0|_p)= (\delta
    A'_0|_p,\delta A'_1|_p)\rho_0(p)+(A'_0|_p,A'_1|_p)\rho_1(p).
     \end{array}\right. 
      \end{equation}
Since $A'$ and $\tilde{A}$ agree at $p$, we then have
\begin{equation}\label{3.7}
 \rho_0=\mathrm{Id}_{2\times 2}.
   \end{equation}
Now, the structure equations of $G$ imply
\begin{equation}\label{3.8}
 dA'_0=\alpha_0^JA'_J,\, dA'_1=\alpha^J_1A'_J,\,
  d\tilde{A}_0=\tilde{\alpha}^J_0\tilde{A}_J,\,
   d\tilde{A}_1=\tilde{\alpha}^J_1\tilde{A}_J,\,
    J=0,\dots,5.
\end{equation}
Substituting (\ref{3.8}) into (\ref{3.6}) yields
\[ 
 \tilde{\alpha}^2_0|_p=\alpha^2_0|_p ,\quad
  \tilde{\alpha}^3_0|_p=\alpha^3_0|_p,\quad
   \tilde{\alpha}^4_0|_p=\alpha^4_0|_p,\quad
    \tilde{\alpha}^2_1|_p=\alpha^2_1|_p,\quad
     \tilde{\alpha}^3_1|_p=\alpha^3_1|_p  
\]
 and
\[ 
 \rho_1(p) = \left( \begin{array}{cc}
  \tilde{\alpha}^0_0|_p-\alpha^0_0|_p &  \tilde{\alpha}^0_1|_p-\alpha^0_1|_p \\
   \tilde{\alpha}^1_0|_p-\alpha^1_0|_p &  \tilde{\alpha}^1_1|_p-\alpha^1_1|_p
    \end{array}\right).
\]
 Since $p$ has been chosen arbitrarily, we can
conclude that the equations
\[
 \tilde{\alpha}^2_0=\alpha^2_0 ,\quad
  \tilde{\alpha}^3_0=\alpha^3_0,\quad
   \tilde{\alpha}^4_0=\alpha^4_0,\quad
    \tilde{\alpha}^2_1=\alpha^2_1,\quad \tilde{\alpha}^3_1=\alpha^3_1
\]
are identically satisfied on $U$. This gives the required result.
\end{proof}

\bigskip
As an application of Theorem \ref{teorema2} we have

\begin{cor} 
If $f,\tilde{f} : M^2 \to \Lambda$ are first order deformations of each other
 then $f$ is a Legendre immersion if and only if $\tilde{f}$ is a Legendre
immersion.
\end{cor}
\bigskip
\subsubsection{Second order deformations} 

We begin by proving the following

\begin{thm}\label{teorema3}
Let $f,\tilde{f}:M\to \Lambda$ be two nondegenerate
Legendre immersions. Then $f$ and $\tilde{f}$ are second order
deformations of each other if and only if there exists a bundle
isomorphism
 \[
  \mathcal{B}:\mathcal{F}(f)\to \mathcal{F}(\tilde{f})
   \]
 such that
\begin{equation}\label{3.10}
\tilde{\omega}^3_0=\mathcal{B}^*(\omega^3_0),\quad
\tilde{\omega}^2_1=\mathcal{B}^*(\omega^2_1).
\end{equation}
\end{thm}

\begin{proof}
Let $U\subset M$ be any coordinate neighborhood of $M$ and let
$A:U\to G$ be a canonical frame along $f$. We show that
a first order deformation $B : M^2 \to G$ is of the second order if and
only if $\tilde{A} = BA : U \to G$ is a canonical frame along $\tilde{f}$ 
such that
\[ 
     \tilde{\alpha}^3_0=\alpha^3_0,\quad 
      \tilde{\alpha}^2_1=\alpha^2_1.
       \]
From Theorem \ref{teorema2}, we know that $\tilde{A}$ is such that
\begin{equation}\label{3.11} 
 \tilde{\alpha}^3_0=\alpha^3_0,\quad
  \tilde{\alpha}^2_1=\alpha^2_1.
    \end{equation}
We also know that the frame fields $\tilde{A}$ and $A'=B(p)A$ satisfy
\begin{equation}\label{3.12}
 \left\{ \begin{array}{ll}
  (\tilde{A}_0,\tilde{A}_1)|_p=(A'_0,A'_1),\\ (\delta \tilde{A}_0,
   \delta \tilde{A}_0)|p= (\delta A'_0,\delta
    A'_1)|_p+(A'_0,A'_1)|_p\rho_1
     \end{array}\right. 
      \end{equation}
where
\[ 
 \rho_1 = \left( \begin{array}{cc}
  (\tilde{\alpha}^0_0-\alpha^0_0)|_p &  (\tilde{\alpha}^0_1-\alpha^0_1)|_p \\
   ( \tilde{\alpha}^1_0-\alpha^1_0)|_p &
    (\tilde{\alpha}^1_1-\alpha^1_1)|_p
     \end{array}\right).
      \]
Lemma \ref{lemma1} implies that $B$ is a second
order deformation if and only if, for every $p\in U$, there exists
\[ 
 \rho_2 = \left( \begin{array}{cc} \sigma^0_0 & \sigma^0_1\\
  \sigma^1_0&\sigma^1_1
   \end{array}\right)\in S^2(M)|_p\otimes \gl(2,\R)
\]
such that
\begin{equation}\label{3.13}
 (\delta^2 \tilde{A}_0,\delta^2 \tilde{A}_1)|_p=(\delta^2 A'_0,
  \delta^2 A'_1)|_p+2(\delta A'_0,\delta A'_1)|_p
   \rho_1+(A'_0,A'_1)|_p\rho_2.
    \end{equation}
Equation (\ref{3.13}), when written out, gives
\begin{equation} \label{3.14}
 \left\{
  \begin{array}{ll} 
 \begin{aligned}
  \delta^2\tilde{A}_0|_p &= \delta^2A'_0|_p+2\delta
   A'_0|_p(\tilde{\alpha}^0_0-\alpha^0_0)|_p+\delta
    A'_1|_p(\tilde{\alpha}^1_0-\alpha^1_0)|_p \\ &\quad +A'_0|_p\sigma^0_0+A'_1|_p\sigma^1_0,
    \end{aligned}\\
     \begin{aligned}
      \delta^2\tilde{A}_1|_p &=\delta^2A'_1|_p+2\delta
       A'_0|_p(\tilde{\alpha}^0_1-\alpha^0_1)|_p+\delta
        A'_1|_p(\tilde{\alpha}^1_1-\alpha^1_1)|_p \\ &\quad+A'_0|_p\sigma^0_1+A'_1|_p\sigma^1_1.
        \end{aligned}
         \end{array}
          \right.
           \end{equation}
On the other hand, from the Maurer-Cartan equations of $\tilde{A}$ and the fact that $\tilde{A}|_p=A'|_p$, we
compute
\begin{equation}\label{3.15}
 \left\{ \begin{array}{ll}
  \delta^2\tilde{A}_0|_p=\gamma^0_0|_pA'_0|_p+\gamma_0^1|_pA'_1|_p+\gamma_0^2|_pA'_2|_p+\gamma_0^3|_pA'_3|_p+
   \gamma^5_0|_pA'_5|_p,\\
    \delta^2\tilde{A}_1|_p=\gamma^0_1|_pA'_0|_p+\gamma_1^1|_pA'_1|_p+\gamma_1^2|_pA'_2|_p+\gamma_1^3|_pA'_3|_p+
     \gamma^4_1|_pA'_4|_p,
      \end{array}\right. 
       \end{equation}
where
\[
\left\{ \begin{array}{rcl}
\gamma_0^0 &=& \delta\tilde{\alpha}^0_0+\tilde{\alpha}^0_0\tilde{\alpha}^0_0+\tilde{\alpha}^1_0\tilde{\alpha}^0_1
+\tilde{\alpha}^3_0\tilde{\alpha}^0_3,\\
\gamma_0^1 &=& \delta\tilde{\alpha}^1_0+\tilde{\alpha}^0_0\tilde{\alpha}^1_0+\tilde{\alpha}^1_0\tilde{\alpha}^0_1+
\tilde{\alpha}^3_0\tilde{\alpha}^1_3,\\
\gamma_0^2 &=& \tilde{\alpha}^1_0\tilde{\alpha}^2_1+\tilde{\alpha}^3_0\tilde{\alpha}^2_3,\\
\gamma_0^3 &=& \delta
\tilde{\alpha}^3_0+\tilde{\alpha}^0_0\tilde{\alpha}^3_0,\\
\gamma_0^5 &=& \tilde{\alpha}^3_0\tilde{\alpha}^3_0,\\
\gamma_1^0 &=& \delta\tilde{\alpha}^0_1+\tilde{\alpha}_1^0\tilde{\alpha}_0^0+\tilde{\alpha}_1^1\tilde{\alpha}_1^0+
  \tilde{\alpha}_1^2\tilde{\alpha}_2^0,\\ 
\gamma_1^1 &=& \delta
 \tilde{\alpha}^1_1+\tilde{\alpha}_1^0\tilde{\alpha}^1_0+\tilde{\alpha}^1_1\tilde{\alpha}^1_1+
  \tilde{\alpha}^2_1\tilde{\alpha}^1_2,\\ 
\gamma_1^2 &=& \delta
 \tilde{\alpha}^2_1 + \tilde{\alpha}^1_1\tilde{\alpha}^2_1,\\
\gamma_1^3 &=& \tilde{\alpha}^0_1\tilde{\alpha}^3_0+\tilde{\alpha}^2_1\tilde{\alpha}^3_2,\\
\gamma_1^4 &=& \tilde{\alpha}^2_1\tilde{\alpha}^2_1.
\end{array}\right. 
\]
Using the Maurer-Cartan equations of $A'$, (\ref{3.14}) can be written as
\begin{equation}\label{3.16}
 \left\{ \begin{array}{ll}
  \delta^2\tilde{A}'_0|_p=\beta^0_0|_pA'_0|_p+\beta_0^1|_pA'_1|_p+\beta_0^2|_pA'_2|_p+\beta_0^3|_pA'_3|_p+
   \gamma^5_0|_pA'_5|_p,\\
    \delta^2\tilde{A}'_1|_p=\beta^0_1|_pA'_0|_p+\beta_1^1|_pA'_1|_p+\beta_1^2|_pA'_2|_p+\beta_1^3|_pA'_3|_p+
     \beta^4_1|_pA'_4|_p,
      \end{array}\right. 
       \end{equation}
where
\[ 
 \left\{ \begin{array}{rcl}
  \beta_0^0 &=& (\delta\alpha^0_0+\alpha^0_0\alpha^0_0+\alpha^1_0\alpha^0_1+\alpha^3_0\alpha^0_3+
   2\alpha^0_0(\tilde{\alpha}^0_0-\alpha^0_0)+2\alpha^0_1(\tilde{\alpha}^1_0-\alpha^1_0)+\sigma^0_0)|_p,\\
    \beta_0^1 &=& (\delta \alpha^1_0+ \alpha^0_0 \alpha^1_1 +
     \alpha^3_0\alpha^1_3+2\alpha^1_0(\tilde{\alpha}^0_0-\alpha^0_0)+
      2\alpha^1_1(\tilde{\alpha}^1_0-\alpha^1_0)+\sigma^1_0)|_p,\\
\beta_0^2 &=& (\alpha^1_0
  \alpha^2_1+2\alpha^2_1(\tilde{\alpha}^1_0-\alpha^1_0))|_p,\\
\beta_0^3 &=& (\delta
  \alpha^3_0+\alpha^0_0\alpha^3_0+2\alpha^3_0(\tilde{\alpha}^0_0-\alpha^0_0))|_p,\\
\beta_0^5 &=& (\alpha^3_0\alpha^3_0)|_p,\\ 
\beta_1^0 &=& (\delta
  \alpha^0_1+\alpha^0_1\alpha^0_0+\alpha^1_1\alpha^0_1+\alpha^2_1\alpha^0_2
   +2\alpha^0_0(\tilde{\alpha}^0_1-\alpha^0_1)+
   2\alpha^0_1(\tilde{\alpha}^1_1-\alpha^1_1)+\sigma^0_1)|_p,\\
\beta^1_1 &=& (\delta
  \alpha^1_1+\alpha^0_1\alpha^1_0+\alpha^1_1\alpha^1_1+\alpha^2_1\alpha^1_2+2\alpha^1_0(\tilde{\alpha}^0_1-\alpha^0_1)+
   2\alpha^1_1(\tilde{\alpha}^1_1-\alpha^1_1)+\sigma^1_1)|_p,\\
\beta_1^2 &=& (\delta
  \alpha^2_1+\alpha^1_1\alpha^2_1+2\alpha^2_1(\tilde{\alpha}^1_1-\alpha^1_1))|_p,\\
\beta_1^3 &=& (\alpha^0_1\alpha^3_0+2\alpha^3_0(\tilde{\alpha}^0_1-\alpha^0_1))|_p,\\
\beta_1^4 &=& (\alpha^2_1\alpha^2_1)|_p.
\end{array}\right.
\]
Form (\ref{3.15}), (\ref{3.16}) and
(\ref{3.11}), it follows that $B$ is a second order deformation if and only if
\begin{equation}\label{3.17}
 \left\{ \begin{array}{rcl}
  \alpha^2(\tilde{\alpha}^1_0-\alpha^1_0) &=& -\alpha^1\tilde{\alpha}^3_2,\\
   \alpha^1(\tilde{\alpha}^0_1-\alpha^0_1) &=& \alpha^2\tilde{\alpha}^3_2,\\
    \alpha^2(\tilde{\alpha}^1_1-\alpha^1_1)&=& 0,\\
     \alpha^1(\tilde{\alpha}^0_0-\alpha^0_0) &= &0.
      \end{array}\right. 
       \end{equation}
The last two equations of (\ref{3.17}) yield
\begin{equation} \label{3.18}
 \tilde{\alpha}^0_0=\alpha^0_0,\quad
  \tilde{\alpha}^1_1=\alpha^1_1.
   \end{equation}
Computing the exterior derivatives of (\ref{3.18})
and using the structure equations, we obtain
\begin{equation}\label{3.19}
 \tilde{\alpha}^1_0\wedge \alpha^2_1 = \tilde{\alpha}^0_1\wedge
  \alpha^3_0=0.
   \end{equation}
Differentiating the equations in (\ref{3.17}), and using again the structure equations, we
see that
\[ 
 \left\{ \begin{array}{ll} \alpha^2\wedge
  (\tilde{\alpha}^1_0-\alpha^1_0)+\alpha^1\wedge
   \tilde{\alpha}^3_2=0,\\ \alpha^1 \wedge
    (\tilde{\alpha}^0_1-\alpha^0_1)-\alpha^2\wedge
     \tilde{\alpha}^3_2=0.
      \end{array}\right.
       \]
This implies that $\tilde{\alpha}^3_2$ vanishes identically, 
Furthermore, from the first two equations of (\ref{3.17}) we get
\[
 \tilde{\alpha}^1_0-\alpha^1_0=\tilde{\alpha}^0_1-\alpha^0_1=0.
  \]
This yields the required result.

\end{proof}

Taking into account \eqref{2.3.4}, we have

\begin{cor}\label{corollario2}
Let $f,\tilde{f}$ be nondegenerate Legendre surfaces and
let $A$ be a canonical frame along $f$. Then $f,\tilde{f}$ are
non-trivial second order deformations of each other if and only if
there exists a normal frame $\widetilde{A}$ along $\widetilde{f}$
such that
\begin{gather}
 \alpha^1 = \widetilde{\alpha}^1,\quad 
   \alpha^2  =  \widetilde{\alpha}^2, \\ 
   \label{3.1.2}
    \alpha^0_3-\widetilde{\alpha}^0_3=w_{1}\alpha^1, \quad
     \alpha^1_2-\widetilde{\alpha}^1_2=w_{2}\alpha^2,\quad
      \alpha^0_4-\widetilde{\alpha}^0_4=-w_{2}\alpha^1+w_{1}\alpha^2,
       \end{gather} 
for smooth functions $w_{1}$, $w_{2}$ such that $(w_1)^2 + (w_2)^2\neq 0$.
\end{cor}

\begin{remark}
If $f$, $\tilde{f}$ are non-trivial second order deformations of each other, then 
from the structure equations of the canonical frames $A$ and $\widetilde{A}$, it follows 
that the $\g$-valued
$1$-form $\eta=\widetilde{\alpha} - \alpha$ satisfies
\begin{equation}\label{d-eta}
  \boxed{d\eta+\alpha\wedge \eta + \eta\wedge \alpha =0.}
   \end{equation} 
\end{remark}

\medskip

We can summarize the previous results in the following.
 
\begin{thm}\label{main-theorem}
 Let $f,\tilde{f} : M^2 \to \Lambda$ be non-congruent, nondegenerate Legendre immersions. 
Then, the following statements are equivalent:

\begin{enumerate}

\item $f$ and $\tilde{f}$ are non-trivial second order deformations of each other.

\item there exists a bundle isomorphism $\mathcal{B}:\mathcal{F}(f)\to
     \mathcal{F}(\tilde{f})$ such that
\[      
 \omega^1=\mathcal{B}^*(\tilde{\omega}^1), \quad \omega^2=\mathcal{B}^*(\tilde{\omega}^2).
  \]
 
\item $f$ and $\tilde{f}$ induce the same canonical coframe on $M$.

\item $f$ and $\tilde{f}$ have the same quotient of cubic to quadratic forms, that is,
\begin{equation}\label{cubic-to-quadratic}
 \Psi/\Phi = \widetilde{\Psi}/\widetilde{\Phi}.
  \end{equation}

\end{enumerate}
\end{thm}

\medskip
\noindent The equivalence of (1), (2), and (3) is a direct consequence of Theorem 
\ref{teorema3}. As for the equivalence with (4), if \eqref{cubic-to-quadratic} holds,
there exist canonical frames $A$ and $\widetilde{A}$ along $f$ and $\widetilde{f}$, respectively, 
such that $\widetilde{\alpha}^1=\alpha^1$, $\widetilde{\alpha}^2=\alpha^2$. 
Since $f,\widetilde{f}$ are not congruent,
this yields \eqref{3.1.2}.

\medskip

\begin{ex}[Isothermic nets]\label{ex: iso}
Let $U\subset \R^2$ be a simply connected open set with coordinates $(u,v)$.
A \textit{net} is a two-parameter smooth immersion $F : U\to \R^3$ 
satisfying $F_{uv}\cdot F_u \times  F_v =0$.
A net $F$ is \textit{isothermic} if $F$ is a
conformal map, that is, $(u,v)$ are both principal and isothermal coordinates. 
Isothermic nets parametrize isothermic surfaces. Examples of isothermic
surfaces include quadrics, surfaces of revolution, cones,
cylinders and constant mean curvature surfaces.
The main
local differential invariant of an isothermic net is the
\textit{Calapso potential}, that is the positive function $\varphi$
defined by 
\[
 \varphi^2(du^2+dv^2)=\frac{1}{4}(k_1-k_2)^2 dF\cdot
  dF,
  \]
where $k_1$ and $k_2$ are the principal curvatures. The
Gauss-Codazzi equations imply that $\varphi$ is a solution of the
\textit{Calapso-Rothe} equation:
\[
 \Delta(\varphi^{-1}\varphi_{uv})+2(\varphi^{2})uv=0.
  \]
Denote by
$f:U\to \Lambda$ the Legendre lift of $F$ and
let $\psi=\log(\varphi)$. If $\psi_u\psi_v\neq 0$, then $f$ is nondegenerate and the
corresponding canonical coframe takes the form 
\[
 \alpha^1=\sqrt[3]{ \psi_u (\psi_v)^2}\,dv, \quad
  \alpha^2=\sqrt[3]{(\psi_u)^2 \psi_v}\,du.
   \] 
Let $W$ be a smooth function such that
$\alpha_{\varphi}=d\left(e^{2\psi}W\right)$, where $\alpha_{\varphi}$ is
the closed $1$-form
\[
\begin{split}
 \alpha_{\varphi} &= -e^{2\psi}\left(\frac{1}{2}\left(e^{-2\psi}\Delta
  \psi \right)_u+2\psi_u \left(1+e^{-2\psi} \Delta \psi \right
   )\right)du \\ & \quad  + e^{2\psi}\left(\frac{1}{2}\left(e^{-2\psi}\Delta
    \psi\right)_v+2\psi_v\left(1+e^{-2\psi}\Delta
     \psi\right)\right)dv.
      \end{split}
       \]
The Calapso potential and the function $W$ give a complete set of invariants for the
isothermic net with respect to the conformal group \cite{Be, MN5}, that is, $\varphi$
and $W$ determine $F$ up to a conformal transformation.
One of the most important features of isothermic nets is the existence of a \textit{spectral transformation}.
This transformation was independently discovered by P. Calapso and L. Bianchi, who introduced it as the
T-\textit{transformation}. Given a real constant $m\in \R$, the $T_m$-transform $T_m(F)$ of $F$ can be 
characterized, up to
conformal transformations, by
\[
 \varphi_{T_m(F)}=\varphi_F,\quad W_{T_m(F)}=W_F+me^{-2\psi}.
  \]
Thus the Legendre lifts of the $T$-transforms of $F$ have the same canonical coframe and are not
congruent (see Remark \ref{rk: isothermic} and \eqref{sezione-parallela}). This shows that 
the Legendre lifts of isothermic nets are deformable.
\end{ex}
\medskip

\begin{ex}[$L$-Isothermic nets]\label{ex: L-iso}
Another class of deformable
surfaces is given by the Legendre lifts of $L$-isothermic nets. A
net $F:U\to \R^3$ is said to be \textit{$L$-isothermic} if 
the Gauss map ${n}:U\to S^2$ is conformal with respect to the 
third fundamental form, that is, $(u,v)$ are 
principal coordinates which are isothermal with respect to the
third fundamental form.
$L$-isothermic nets parametrize the class of $L$-\textit{isothermic surfaces}. 
Examples include minimal surfaces
in $\R^3$ and molding surfaces \cite{MN3}. The study of $L$-isothermic surfaces
goes back to the work of Blaschke and presents many analogies
with that of isothermic surfaces. For instance, $L$-isothermic
surfaces admit a spectral transformation which is the analogue of
the T-transformation for isothermic surfaces \cite{MN1}. We briefly recall some
basic properties of $L$-isothermic nets and show that
their Legendre lifts are deformable. The
\textit{Blaschke potential} of $F$ is the positive function
$\varphi$ defined by
\[
 \varphi^2(du^2+dv^2)=\frac{1}{4}\left(\frac{1}{k_1}-\frac{1}{k_2}\right)^2
  d \,{n} \cdot d\, {n}.
   \] 
In this case the compatibility condition arising from the Gauss-Codazzi equations is the 
\textit{Bla\-schke equation}
\begin{equation}\label{3.22}
 \Delta\left(\frac{1}{\varphi}\varphi_{uv}\right)=0.
  \end{equation}
Let $f$ be the Legendre lift of $F$ and $\psi=\log(\varphi)$.
In the nondegenerate case, which amounts to $\psi_u\psi_v\neq 0$, the canonical coframe of $f$ can 
be written as 
\[
 \alpha^1=\sqrt[3]{(\psi_u)^2 \psi_v}du, \quad    
  \alpha^2=\sqrt[3]{\psi_u (\psi_v)^2}dv.
   \]
In addition to the Blaschke potential, the other local differential invariant of $F$ 
is a smooth function $W$ defined by $\alpha_{\varphi} = d\left(e^{2\psi}W\right)$,
where 
\[
 \begin{split}
  \alpha_{\varphi} &= -e^{2\psi}\left(\frac{1}{2}\left(e^{-2\psi}\Delta
   \psi \right)_u+2\psi_u \left(e^{-2\psi} \Delta \psi \right
    )\right)du \\ & \quad  + e^{2\psi}\left(\frac{1}{2}\left(e^{-2\psi}\Delta
     \psi\right)_v+2\psi_v\left(e^{-2\psi}\Delta
      \psi\right)\right)dv.
       \end{split}
        \]
Given $m\in \R$, the $T_m$-transform $T_m(F)$ of $F$ can be
characterized (up to Laguerre contact transformations) by
\[
 \varphi_{T_m(F)}=\varphi_F,\quad W_{T_m(F)}=W_F+me^{-2\psi}.
  \] 
The Legendre lifts of the $T$-transforms of $F$ have then the same
canonical coframe and are not congruent (see Remark \ref{rk: isothermic} and 
\eqref{sezione-parallela}), from which follows that the Legendre lift of an
$L$-isothermic net is deformable.
\end{ex}

\bigskip

\section{Infinitesimal deformations and deformable surfaces}\label{s: infinitesimal-defo}

\textit{Infinitesimal deformations.} On the one hand, if $f$ and $\widetilde{f}$ are non-trivial 
deformations of each other, then $\eta = \widetilde{\alpha}-{\alpha}$ never vanishes and
according to \eqref{d-eta}
\begin{equation}\label{infinitesimal-displacement}
 \delta|_U := A\eta A^{-1}
  \end{equation}
is a closed 1-form with values in $\g$, for every canonical frame $A$ along $f$.
Moreover,
\begin{equation}\label{displacement}
  D|_U : = [\widetilde{A}A^{-1}]
  \end{equation}
defines a smooth map $D : M^2 \to G/ \Z_2$ such that
\[
 \quad D^{-1}dD=\delta.
   \] 

 On the other hand, let $f : M^2 \to \Lambda$ be a nondegenerate Legendre immersion and define
$\eta(w_1,w_2)\in \Omega^1(U)\otimes \g$ by
\begin{equation}\label{etaw1w2}
 \eta(w_1,w_2) = \left(\!
   \begin{array}{cccccc}
  0 & 0 & 0 & w_1\alpha^1& -w_2\alpha^1 + w_1 \alpha^2 &0 \\
  0 & 0 & w_2\alpha^2 &0&0&w_2\alpha^1 - w_1\alpha^2\\ 
  0 & 0 &  0          &0&w_2 \alpha^2 &0 \\ 
  0&0&0&0&0&w_1\alpha^1\\
  0&0&0&0&0&0\\ 
  0&0&0&0&0&0
  \end{array}
   \right),
   \end{equation}
for smooth functions $w_1,w_2$. Note that $\eta$ takes values in the abelian subalgebra
\[
 \a = \{ T\in \g \, |\,
   T(\epsilon_0)=T(\epsilon_1)=0,\, T(\epsilon_2)\varpropto
   \epsilon_1,\, T(\epsilon_3)\varpropto \epsilon_0 \}.
    \]
From \eqref{d-eta}, it follows that the $1$-form $A \eta A^{-1}\in \Omega^1(U)\otimes \g$ is
independent of $A$. Thus, there exists $\delta\in \Omega^1(M)\otimes \g$ such 
that $\delta |_U=A \eta A^{-1}$.

\begin{defn}[Infinitesimal deformations]
We say that $\eta$ is an \textit{infinitesimal deformation}
of $f$ if $\delta$ is closed. Let $\Delta_{f}$ denote the
set of infinitesimal deformations of $f$.
\end{defn}
 
We are now in a position to characterize the Legendre surfaces which admit non-trivial 
deformations in terms of infinitesimal deformations.

\begin{thm}
 Let $M^2$ be simply connected. 
  Then a nondegenerate Legendre immersion $f : M^2 \to \Lambda$
   admits non-trivial deformations if and only if $\Delta_{f}\neq 0$. 
\end{thm}

\begin{proof}
If $\widetilde{f}$ is a non-trivial deformation of $f$, then $\eta=\widetilde\alpha -\alpha$ defines
a non-zero
infinitesimal deformation. Conversely, let $\eta$
be a non zero infinitesimal deformation. Then, $\delta\in
\Omega^1(M)\otimes \g$ is a non-zero closed $1$-form
and there exists $D:M\to G$
such that $D^{-1}dD=\delta$. Note that the map $\widetilde{f}: p\in M\mapsto D(p) \cdot f(p)\in
\Lambda$ is a non-trivial deformation of $f$.
\end{proof}

\begin{remark}
Note that for every non-zero $\eta\in \Delta_{f}$ there exists a 
non-trivial deformation
$\widetilde{f}$
which is uniquely defined by $f$ and the corresponding infinitesimal deformation,
 up to the action of $G$.
Moreover,
as $\Delta_{f}$ is a real vector space,
given $\eta\in \Delta_{f}$ and $r\in \R$, $r\eta$ is
another infinitesimal deformation. Therefore, the deformations of
$f$ arise in one-parameter families. In other words, deformable
surfaces do have a \textit{spectral transformation}. This suggests 
the existence of a B\"acklund transformation for deformable surfaces.
\end{remark}

\textit{Infinitesimal deformations and parallel sections.} Let $f : M^2 \to \Lambda$ be a nondegenerate
Legendre immersion and consider the $\gl(3,\R)$-valued $1$-form
\begin{equation} \label{3.3.1}
\sigma = \left(\begin{array}{ccc}
   -2(2q_1\alpha^1-q_2\alpha^2)&0&-\alpha^1\\
     0&-2(q_1\alpha^1-2q_2\alpha^2)&\alpha^2\\
      2(p_2-1)\alpha^2&-2(p_1-1)\alpha^1&-3(q_1\alpha^1-q_2\alpha^2)
       \end{array} \right).
        \end{equation}

\medskip
\begin{defn}[$\sigma$-connection]
The form \eqref{3.3.1} defines a linear connection 
\[
 D^{\sigma}w := dw + \sigma w 
  \]
on the trivial bundle $M^2\times \R^3$, for each smooth function $w : M^2 \to \R^3$. 
$D^{\sigma}$ is referred to as the \textit{$\sigma$-connection} of $f$. 
By $\mathcal{P}_{f}$ we denote the vector space consisting of all \textit{parallel
sections} of the $\sigma$-connection. 
\end{defn}

A simple computation shows
that the curvature form $\Omega^{\sigma}$ of the $\sigma$-connection is given by
\begin{equation}\label{3.3.3}
  \Omega^{\sigma} = \left(\begin{array}{ccc} 0 & 0 & 0\\ 0 & 0 & 0
   \\ 2\partial_1 p_2 & 2\partial_2 p_1& 3(p_2-p_1)
    \end{array} \right)\alpha^1\wedge \alpha^2,
     \end{equation}
where for a smooth function $g : M\to \R$ we write $dg= \partial_1g \,\alpha^1 + \partial_2g\, \alpha^2$.

\medskip
 For every $w = (w_1,w_2,w_3) : M^2 \to \R^3$, let $\eta(w_1,w_2)$ be defined by \eqref{etaw1w2}. 

\begin{prop} 
A nondegenerate Legendre immersion $f$ is deformable if and only if there exists a parallel
section of the $\sigma$-connection. 
 Moreover, the mapping
\[
 w \in \mathcal{P}_{f}\to \eta(w_1,w_2)\in
  \Delta_{f}
   \] 
is an isomorphism of vector spaces.
\end{prop}

\begin{proof} Let $\widetilde{f}$ be a deformation of $f$.
Then $\alpha^1=\widetilde{\alpha}^1$, $\alpha^2=\widetilde{\alpha}^2$ and
\[
 \alpha^0_3-\widetilde{\alpha}^0_3=w_1\alpha^1,\quad
  \alpha^1_2-\widetilde{\alpha}^1_2=w_2\alpha^2,\quad
   \alpha^0_4-\widetilde{\alpha}^0_4=-w_2\alpha^1+w_1\alpha^2,
    \] 
for smooth functions $w_1,w_2$ such that $(w_1)^2+(w_2)^2 \neq 0$. 
Differentiating and using the structure equations, we get
\[
 \begin{array}{lll}
   d(w_1\alpha^1)=-\alpha^0_0\wedge w_1\alpha^1,\\
    d(w_2\alpha^2)=-w_2\alpha^1_1\wedge \alpha^2,\\
     d(-w_2\alpha^1+w_1\alpha^2)=-(\alpha^0_0+\alpha^1_1)\wedge(-w_1\alpha^1+w_1\alpha^2),
      \end{array}
      \]
which implies
\[
 \begin{array}{lll}
  &(dw_1+2w_1\alpha^0_0)\wedge \alpha^1=0,\\ &
   (dw_2+2w_2\alpha^1_1)\wedge \alpha^2=0,\\
    &(dw_1+2w_1\alpha^0_0)\wedge\alpha^2-(dw_2+2w_2\alpha^1_1)\wedge\alpha^1=0,
     \end{array}
      \]
By Cartan's Lemma, there exists a smooth function
$w_3 : M^2 \to \R$ such that 
\[
 dw_1=-2w_1\alpha^0_0+w_3\alpha^1,\quad dw_2=-2w_2\alpha^1_1-w_3\alpha^2,
  \] 
that is,
 \begin{eqnarray*}
  dw_1-2w_1(2q_1\alpha^1-q_2\alpha^2)-w_3\alpha^1=0,\\
  dw_2-2w_2(q_1\alpha^1-2q_2\alpha^2)+w_3\alpha^2=0.
   \end{eqnarray*}
Taking the exterior derivative of these equations yields
\[
 dw_3+2(p_2-1)\alpha^2w_1-2(p_1-1)\alpha^1w_2-3(q_1\alpha^1-q_2\alpha^2)w_3=0,
  \]
which shows that $(w_1,w_2,w_3)$ is a parallel cross section of the $\sigma$-connection.

The converse follows by observing that $\eta(w_1,w_2)$
defines an infinitesimal deformation if $w=(w_1,w_2,w_3)$ is a parallel section.

\end{proof}

\begin{remark}
Nondegenerate deformable Legendre surfaces may be classified in terms of the
dimension of $\Delta_{f}$. For a generic $f$ the space
$\Delta_{f}$ is one-dimensional. It is not too difficult to
show that surfaces with a three-parameter family of deformations
can be generically obtained as deformations of the Legendre lifts
of molding surfaces in $\R^3$. It is not at all clear if there exist
Legendre surfaces with a two-dimensional family of deformations.
But, if they exist then they can be reconstructed from the
integral manifolds of a Pfaffian system with empty complex characteristic variety. Thus, 
this class is either empty, or it depends on a finite number of parameters.
\end{remark}

\bigskip

\section{Examples}\label{s: examples}

\begin{ex}[Generic deformations]
 Let $f$, $\widetilde{f}$
be deformations of each other and let $\eta$ be the
corresponding infinitesimal deformation. According to the notations of Corollary \ref{corollario2}, 
we say that the deformation 
is \textit{generic} if $w_1w_2$ is nowhere
vanishing. In this case, there
exist local parameters $(u,v)$ on $M$ such that the canonical
coframe takes the form 
\[ 
 \alpha^1=\sqrt[3]{(
  \psi_u)^2 \psi_v}du,\quad \alpha^2=\sqrt[3]{\psi_u
   (\psi_v)^2}dv,
    \] 
where $\psi :M^2 \to \R$ is a smooth
function. 
The $T$-transforms of isothermic and $L$-isothermic nets are examples of generic deformations. 

\begin{remark}\label{rk: isothermic}
In the case of an isothermic, respectively, $L$-isothermic surface,
$e^{\psi}$ is the Calapso, respectively, the Blaschke potential. This can be
seen by applying the reduction procedure recalled in Section \ref{s: basic} 
to their respective conformal and Laguerre canonical frames. 
For this we need to assume the nondegeneracy condition $\psi_u\neq 0$, $\psi_v\neq 0$.
\end{remark}

\noindent A direct computation shows that the $\sigma$-connection
$\sigma=(\sigma_a^b)$ is given by
\[
 \sigma^1_2=\sigma^2_1=0,\quad
  \sigma^1_3=-\alpha^1=-\sqrt[3]{(
   \psi_u)^2 \psi_v}du,\quad
    \sigma^2_3=\alpha^2=\sqrt[3]{\psi_u (\psi_v)^2}dv,
     \]
and
\[
 \left\{\begin{array}{l}
  \sigma^1_1 =
   \frac{4}{3}\left(\frac{\psi_{uu}}{
    \psi_u}+\frac{2 \psi_{uv}}{
     \psi_v}\right)du+\frac{2}{3}\left(\frac{2\psi_{uv}}{\psi_u}+
      \frac{\psi_{vv}}{
       \psi_v}\right)dv\\
       \sigma^2_2= \frac{2}{3}\left(\frac{\psi_{uu}}{\psi_u}+
       \frac{2\psi_{uv}}{
       \psi_v}\right)du + \frac{4}{3}\left(\frac{2\psi_{uv}}{\psi_u}+
      \frac{\psi_{vv}}{\psi_v}\right)dv,\\
  \sigma^3_1=   -\frac{2}{\sqrt[3]{(
\psi_u)^2\psi_v}}\left(\frac{\psi_{uvv}}{\psi_{v}}
  -\frac{\psi_{vv}
  \psi_{uv}}{(\psi_v)^2}\right)dv,\\
  \sigma^3_2= \frac{2}{\sqrt[3]{\psi_u (\psi_v)^2}}\left(\frac{\psi_{uuv}}{\psi_u}
  -\frac{\psi_{uu} \psi_{uv}}{(\psi_u)^2}\right)du
  ,\\
    \sigma^3_3=   \left(\frac{\psi_{uu}}{ \psi_u}+\frac{2\psi_{uv}}{
     \psi_v}\right)du+\left(\frac{2\psi_{uv}}{\psi_u}+\frac{\psi_{vv}}{
      \psi_v}\right)dv
       \end{array}\right. 
        \]
It is a computational matter to verify that the parallel section
$w\in \mathcal{P}_{f}$ associated to $\eta$ is given
by
\begin{equation}\label{sezione-parallela}
  w={^t\left(\sqrt[3]{( \psi_u)^{-4}(
   \psi_v)^{-2}},-\sqrt[3]{( \psi_u)^{-2}(
    \psi_v)^{-4}},2(\psi_u)^{-2}(
     \psi_v)^{-2} \psi_{uv} \right).}
      \end{equation}

\begin{remark}
Note that $w^1, w^2\neq 0$, which 
characterizes such surfaces.
Moreover, in the case of isothermic and $L$-isothermic nets, $w$ originates the one-parameter family of
non-trivial deformations considered in Examples \ref{ex: iso} and \ref{ex: L-iso}.
\end{remark}

\end{ex}

\begin{ex}[Special deformations]
A deformation is said to be \textit{special} if $w_1w_2$ vanishes identically. Deformable surfaces which 
admit a special deformation play the role of the $R_0$ surfaces in projective differential
geometry. Let $w \in \mathcal{P}_{f}$ be the parallel section
associated to the deformation. Then $w^1w^2$ and $w^3$ vanish
identically. Two cases may occur: either $w^1=0$, or else $w^2=0$.
Without loss of generality, we assume that $w^2=0$. From the
structure equations of the canonical frame, it follows that
$f$ admits a special deformation with $w^2=0$ if and only if
$p_2=1$. The degree of generality of this class of Legendre immersions 
will be clear in the last section; we will see that
they are rather special. Notice that on $M^2$ there
exist local parameters $(u,v)$ with respect to which the canonical
coframe is given by 
\[
 \alpha^1=\psi^{2/3}du,\quad
  \alpha^2=\psi^{1/3}dv,
   \]
for a smooth function $\psi$ such
that $w^1=\pm \, \psi^{-2/3}$.
\end{ex}

\begin{ex}[Legendre surfaces with $3$-parameter families of
deformations]\label{esempio4}

In this example we consider the Legendre immersions with flat $\sigma$-connection. 
This example has been discussed by Ferapontov in \cite{Fe3} (see also Finikov \cite{Fi}). 
From (\ref{3.3.3}) we see that
$\Omega^{\sigma}=0$ if, and only if, $p_1=p_2=c$, for a constant 
$c$.
According to Remark \ref{rk: 3-web}, we have

\begin{prop}
The $3$-web defined by the canonical coframe is flat if and only if $p_1=p_2$.
\end{prop}

\begin{remark} 
From this we infer that deformable diagonally cyclidic
surfaces are characterized by having $p_1=p_2=\mbox{const}$. 
\end{remark}

\noindent Since the web-connection is flat, then there exist local
coordinates $(u,v)$ such that
\begin{equation}\label{eq2esempio4}
 \alpha^1=e^{\psi}du,\quad
  \alpha^2=e^{\psi}dv,
\end{equation}
 where $\psi$ is a smooth function. From this we see that
\begin{equation}
  q_1=-\psi_ue^{-\psi},\quad
    q_2=\psi_ve^{-\psi}.
      \end{equation}
This implies
\begin{equation}\label{eq3esempio4}
  \left\{ \begin{array}{ll}
    \alpha^0_0=2\psi_udu+\psi_vdv,\\ \alpha^1_1=\psi_udu+2\psi_vdv
     \end{array}\right. 
      \end{equation}
From the structure equation we deduce that $\psi$ is a
solution of the Liouville equation
\begin{equation}\label{eq4esempio4}
 \psi_{uv}=(1-c)e^{2\psi}.
   \end{equation}
The other compatibility conditions arising from the structure equations are:
\begin{equation}\label{eq5esempio4}
 \left\{ \begin{array}{lll}
  d\alpha^1_2=-\alpha^1_1\wedge \alpha^1_2,\\
   d\alpha^0_3=-\alpha^0_0\wedge \alpha^0_3,\\
    d\alpha^0_4=-(\alpha^0_0+\alpha^1_1)\wedge \alpha^0_4,
      \end{array}\right. 
       \end{equation}
where
\begin{equation}\label{eq6esempio4}
 \left\{ \begin{array}{lll}
  \alpha^0_3=Adu+ce^{\psi}dv,\\ \alpha^1_2=ce_{\psi}du+Bdv,\\
   \alpha^0_4=-Bdu+Adv,
    \end{array}\right. 
     \end{equation}
for two suitable smooth functions $A,B$
(essentially the invariant $r_1$ and $r_2$). It is now a
computational matter to verify that (\ref{eq5esempio4}) and
(\ref{eq6esempio4}) can be written as follows
\begin{equation}\label{eq7esempio4}
 \left\{ \begin{array}{lll}
   A_v=-A\psi_v+3ce^{\psi}\psi_u,\\ B_u=-B\psi_u+3ce^{\psi}\psi_v,\\
    A_u+B_v=-3A\psi_u-3B\psi_v.
     \end{array}\right. 
      \end{equation}
We may rewrite (\ref{eq7esempio4}) in the form
\begin{equation}\label{eq8esempio4}
 \left\{ \begin{array}{ll}
   dA=(R-3A\psi_u)du+(3ce^{\psi}\psi_u-A\psi_u)dv,\\
    dB=(3ce^{\psi}\psi_v-B\psi_u)du-(R+3B\psi_v)dv,
     \end{array}\right. 
      \end{equation}
where $R$ is a suitable smooth function. Differentiating (\ref{eq8esempio4}), we get
\begin{equation}\label{eq9esempio4}
 \left\{ \begin{array}{ll}
   R_u=-R\psi_u-2(1-c)e^{2\psi}B-3ce^{\psi}(\psi_{vv}+4\psi_v^2,\\
    R_v=-R\psi_v+2(1-c)e^{2\psi}A+3ce^{\psi}(\psi_{uu}+4\psi_u^2).
     \end{array}\right. 
      \end{equation}
The compatibility conditions of this system imply 
\begin{equation}\label{eq10esempio4}
 3ce^{\psi}\left(\psi_{uuu}+\psi_{vvv}+10\psi_u\psi_{uu}+10\psi_v\psi_{vv}
  +8(\psi_u^2+\psi_v^2)\right)=0.  
   \end{equation}
Two cases may occur: either $c=0$, or $c\neq0$. In the
first case the only compatibility condition is the Liouville
equation, which may be viewed as the soliton equation of this class
of surfaces: its solutions depend on two arbitrary functions in
one variable. In fact, the general solutions of the Liouville
equation are of the form
\begin{equation}\label{eq11esempio4}
 \left\{ \begin{array}{ll}
  \Phi^{2}=\frac{1}{1-c}\frac{\lambda'\dot{\mu}}{(\lambda+\mu)^2},\quad
   \mathrm{if}\quad c\neq 1,\\ \Phi^2=\lambda'\dot{\mu},\quad
    \mathrm{if}\quad  c=1,
     \end{array}\right. 
       \end{equation}
where $\Phi=e^{\psi}$, $\lambda$ is a function of the variable $u$, and $\mu$ is a function 
of the variable $v$.
Thus, if $c\neq 0$, it follows that $\psi$ is a solution of the
overdetermined system
\begin{equation}\label{eq12esempio4}
 \left\{ \begin{array}{ll}
  \psi_{uuu}+\psi_{vvv}+10\psi_u\psi_{uu}+10\psi_v\psi_{vv}+8(\psi_u^2+\psi_v^2)=0,\\
   \psi_{uv}=(1-c)e^{2\psi}.
    \end{array}\right.
     \end{equation}
If we use the potential $\Phi=e^{\psi}$, then (\ref{eq12esempio4}) means that $\Phi^2$ 
is a function of the form (\ref{eq11esempio4}) such that
\begin{equation} \label{eq13esempio4}
 \left(\frac{1}{\Phi^2}\left(\Phi^2\left(\Phi^2\right)_u\right)_u\right)_u+
  \left(\frac{1}{\Phi^2}\left(\Phi^2\left(\Phi^2\right)_v\right)_v\right)_v=0.
   \end{equation}
Let first examine the case $c=1$. We take $\lambda$ and $\mu$ as a
new variable and we think of $\lambda'$ and $\dot{\mu}$ as a
functions of $\lambda$ and $\mu$ respectively. With respect to
these new coordinates the equation (\ref{eq13esempio4}) is
equivalent to
\begin{equation}\label{eq14esempio4}
 \partial^{3}_{\lambda\lambda\lambda}(\lambda'^3)+\partial^{3}_{\mu\mu\mu}(\dot{\mu}^3)=0.
  \end{equation}
This implies that $\lambda$ and $\mu$ satisfy the ODE
\begin{equation}\label{eq15esempio4}
 \left\{ \begin{array}{ll}
  \left(\frac{d\lambda}{du}\right)^3=P(\lambda),\\
   \left(\frac{d\mu}{dv}\right)^3=Q(\mu),
    \end{array}\right.
     \end{equation}
where $P$ and $Q$ are polynomials of order $\leq 3$ with
the opposite leading coefficients.
\bigskip
One may proceed in a similar fashion also in the general case. Substituting
\begin{equation}
 \Phi^{2}=\frac{1}{1-c}\frac{\lambda'\dot{\mu}}{(\lambda+\mu)^2}
  \end{equation}
in (\ref{eq14esempio4}), we obtain
\begin{equation}\label{eq16esempio4}
 \begin{split}
  \frac{1}{3}(\lambda+\mu)^3 &\left(\frac{d^3\lambda'^3}{d\lambda^3}+ \frac{d^3
   \dot{\mu}^3}{d\mu^3}\right) - 4(\lambda+\mu)^2\left(\frac{d^2\lambda'^3}{d\lambda^2}+\frac{d^2
    \dot{\mu}^3}{d\mu^2}\right) \\ 
     &\quad
      +20(\lambda+\mu)\left(\frac{d\lambda'^3}{d\lambda}+\frac{d
       \dot{\mu}^3}{d\mu}\right)-40(\lambda'^3+\dot{\mu}^3)=0.
        \end{split}
         \end{equation}
Applying the operator
$\partial^6_{\lambda\lambda\lambda\mu\mu\mu}$ to (\ref{eq16esempio4}), we get
 \begin{equation}\label{eq17esempio4}
  \frac{d^6\lambda'^3}{d\lambda^6}+\frac{d^6\dot{\mu}^3}{d\mu^6}=0.
    \end{equation}
This implies that $\lambda'^3$ and $\dot{\mu}^3$ are
polynomials $P(\lambda)$ and $Q(\mu)$ in $\lambda$ and $\mu$,
respectively, of order $\leq 6$. Such polynomials satisfy
(\ref{eq16esempio4}) if and only if $Q(T)=-P(-T)$. Thus, $\lambda$
and $\mu$ satisfy the ODE
  \begin{equation} \label{eq18esempio4}
   \left\{ \begin{array}{ll}
    \left(\frac{d\lambda}{du}\right)^3=P(\lambda),\\
     \left(\frac{d\mu}{dv}\right)^3=-P(-\mu),
      \end{array}\right.
       \end{equation}
where $P$ is a polynomial of degree $\leq 6$.

\bigskip
\textit{The $3$-parameter family of infinitesimal deformations.} We finish this example by 
discussing the $3$-parameter
family of infinitesimal deformations of such surfaces. The
$\sigma$-connection is given by
  \begin{equation}\label{eq19esempio4}
\sigma=\left(\begin{array}{ccc}
  4\psi_udu+2\psi_vdv & 0 & -e^{\psi}du \\
  0 & \psi_udu+4\psi_vdv & e^{\psi}dv \\
  2(c-1)e^{\psi}dv & -2(c-1)e^{\psi}du & 3d\psi
   \end{array}
     \right),
      \end{equation}
where the potential function is given either by
\begin{equation}\label{eq20esempio4}
e^{2\psi}=\lambda' \dot{\mu},
  \end{equation}
or by
\begin{equation}\label{eq21esempio4}
  e^{2\psi}=\frac{1}{1-c}\frac{\lambda' \dot{\mu}}{(\lambda+\mu)^2}.
   \end{equation}
It is now a computational matter to verify that the
parallel sections of the $\sigma$-connection are given by:

\medskip
$\bullet$ If $c=1$ (i.e., $e^{2\psi}=\lambda' \dot{\mu}$)
\begin{equation}
 s_0 e^{-2\psi}\left(\begin{array}{c}
  \frac{1}{\lambda'} \\
  -\frac{1}{\dot{\mu}} \\
     0
     \end{array}\right)
       +s_1 e^{-2\psi}\left(\begin{array}{c}
  \frac{\lambda}{\lambda'} \\
  -\frac{\mu}{\dot{\mu}} \\
  e^{-\psi}
\end{array}\right)+s_2 e^{-2\psi}\mu \lambda \left(\begin{array}{c}
  \frac{\lambda}{\lambda'} \\
  -\frac{\mu}{\dot{\mu}} \\
  2 e^{-\psi}
     \end{array}\right)
        \end{equation}
where $s_0,s_1$ and $s_3$ are real constants. 

\medskip
$\bullet$ If $c\neq 1$ (i.e.,
$e^{2\psi}=\frac{1}{1-c}\frac{\lambda'
\dot{\mu}}{(\lambda+\mu)^2}$), we have
\begin{equation}
s_0e^{-2\psi}\left(\begin{array}{c}
  \frac{-1}{\lambda'} \\
  \frac{1}{\dot{\mu}} \\
\frac{2e^{-\psi}}{\lambda+\mu}
\end{array}\right)
+s_1e^{-2\psi}\left(\begin{array}{c}
  \frac{\lambda}{\lambda'} \\
  \frac{\mu}{\dot{\mu}} \\
  \frac{e^{-\psi}(\mu-\lambda)}{\lambda+\mu}
\end{array}\right)
+s_2e^{-2\psi}\left(\begin{array}{c}
  \frac{\lambda^2}{\lambda'} \\
  \frac{-\mu^2}{\dot{\mu}} \\
  \frac{2\mu\lambda e^{-\psi}}{\lambda+\mu}
\end{array}\right),
\end{equation}
\noindent where $s_0$,$s_1$ and $s_2$ are real constants.
\end{ex}

\bigskip

\section{The differential system of a deformation}\label{s: diff-system}

It  was shown in \cite{J} that the problems of $k$-th order deformation are equivalent
to solving certain exterior differential systems on appropriate spaces. 
Naturally, for each concrete homogeneous space there is a specific problem to solve. 
We shall derive this result in the case at hand.

\medskip
Let $P=(G/\Z_2)\times \R^6 \times \R^3$ and denote by
$(q_1,q_2,p_1,p_2,r_1,r_2)$ and $(w_1,w_2,w_3)$ the
coordinates on $\R^6$ and $\R^3$, respectively.
Let $(\omega^I_J)$ be the Maurer-Cartan forms on
$G/\Z_2$ and put $\alpha^1=\omega^3_0$,
$\alpha^2=\omega^2_1$. On $P$, we consider the exterior
differential $1$-forms $\eta^1, \dots,\eta^{16}$ defined by
\begin{equation}\label{Generators1}зк
 \left\{ \begin{array}{lcllcllcl} 
  \eta^1 &=&\omega^4_0, &\eta^2 &=& \omega^2_0, &\eta^3 & =&\omega^3_1,\\ 
   \eta^4 &=&\omega^3_2, & \eta^5 &=&\omega^1_0-\alpha^2, & \eta^6 &=& \omega^0_1-\alpha^1,\\ 
    \eta^7 &=&\omega^0_2,& \eta^8 &=&\omega^1_3, &&&
     \end{array}\right. 
      \end{equation}
\begin{equation}\label{Generators2}
 \left\{ \begin{array}{lcl} 
   \eta^9 & = & \omega^0_0+2q_1\alpha^1-q_2\alpha^2,\\ 
    \eta^{10}& = &\omega^1_1+q_1\alpha^1-2q_2\alpha^2,
     \end{array}\right. 
      \end{equation}
\begin{equation}\label{Generators3}
 \left\{ \begin{array}{lcl} 
  \eta^{11}&=&\omega^0_3-r_1\alpha^1-p_2\alpha^2,\\ 
   \eta^{12}&=&\omega^1_2-p_1\alpha^1-r_2\alpha^2,\\ 
    \eta^{13}&=&\omega^0_4+r_2\alpha^1-r_1\alpha^2,
      \end{array}\right. 
      \end{equation}
\begin{equation}\label{Generators4}
 \left\{ \begin{array}{lcl} 
  \eta^{14}&=&dw_1-2w_1(2q_1\alpha^1-q_2\alpha^2)-w_3\alpha^1,\\ 
   \eta^{15}&=&dw_2-2w_2(q_1\alpha^1-2q_2\alpha^2)+w_3\alpha^2,\\ 
    \eta^{16}&=&dw_3+2w_1(p_2-1)\alpha^2-2w_2(p_1-1)\alpha^1-3w_3(q_1\alpha^1-q_2\alpha^2).
     \end{array}\right. 
      \end{equation}

\begin{defn}
Let $(\mathcal{I},\alpha^1\wedge\alpha^2)$ be the exterior
differential system on $P$ generated by $\eta^1,\dots,\eta^{16}$
with the independence condition $\alpha^1\wedge \alpha^2\neq 0$.
We call $(\mathcal{I},\alpha^1\wedge\alpha^2)$ the \textit{differential
system of a deformation}.
\end{defn}

\begin{remark}
The integral manifolds of $(\mathcal{I},\alpha^1\wedge\alpha^2)$ are $2$-dimensional 
immersed surfaces $([A],q,p,r,w):M\to P$ such that

\begin{itemize}

\item $f=[A_0\wedge A_1]:M\to \Lambda$ is a nondegenerate Legendre immersion.

\item $[A]:M\to G/\Z_2$ is the canonical frame along $M$.

\item $q,p, r:M\to \R^2\times \R^2\times \R^2$ are the invariant functions of $f$.

\item $w:M\to \R^3$ is a parallel section of the
$\sigma$-connection of $f$.

\end{itemize}
Thus, the deformations of a
nondegenerate Legendre immersion may be regarded as the integral manifolds of
the differential system $(\mathcal{I},\alpha^1\wedge\alpha^2)$.
\end{remark} 

From the
Maurer--Cartan equations we obtain the \textit{quadratic equations} of the system, 
which are (modulo $\mathcal{I}$)
 \begin{equation} \label{QE1}
  d\alpha^1\equiv -q_2\alpha^1\wedge \alpha^2,\quad
    d\alpha^2\equiv -q_1\alpha^1\wedge \alpha^2, 
    \end{equation}
\begin{equation} 
  d\eta^1\equiv\dots \equiv d\eta^8 \equiv 0,
   \end{equation}
\begin{equation}\label{QE2}
 \left\{ \begin{array}{lcl} 
  d\eta^9 & \equiv & 2dq_1\wedge \alpha^1-dq_2\wedge \alpha^2+(-1+p_2-q_1q_2)\alpha^1\wedge \alpha^2,\\ 
   d\eta^{10} &\equiv &dq_1\wedge \alpha^1-2dq_2\wedge\alpha^2+(1-p_1+q_1q_2)\alpha^1\wedge \alpha^2,
    \end{array}\right.
     \end{equation}
\begin{equation}\label{QE3}
 \left\{ \begin{array}{lcl} 
  d\eta^{11} &\equiv &-dr_1\wedge\alpha^1-dp_2\wedge\alpha^2+(2r_1q_2+3q_1p_2)\alpha^1\wedge \alpha^2 ,\\ 
   d\eta^{12} &\equiv &-dp_1\wedge\alpha^1-dr_2\wedge\alpha^2+(3p_1q_2+2r_2q_1)\alpha^1\wedge\alpha^2,\\ 
    d\eta^{13} &\equiv &dr_2\wedge \alpha^1-dr_1\wedge\alpha^2 + 4(q_1r_1-q_2r_2)\alpha^1\wedge \alpha^2,
     \end{array}\right.
       \end{equation}
\begin{equation}\label{QE4}
  d\eta^{14}\equiv -2w_1d\eta^9,\quad d\eta^{15}\equiv -2w_2d\eta^{10}
    \end{equation}
\begin{equation}\label{QE5}
  d\eta^{16}\equiv -w_3(d\eta^9+d\eta^{10})+2w_1dp_2\wedge
   \alpha^2-2w_2dp_1\wedge\alpha^1+3w_3(p_2-p_1)\alpha^1\wedge\alpha^2.
    \end{equation}
From this, we see that the differential ideal $\mathcal{I}$ is
algebraically generated by the $1$-forms $\eta^1,\dots,\eta^{16}$
and by the differential $2$-forms
\begin{equation}\label{QuadraticGenerators}
  \left\{ 
\begin{array}{lcl} 
\Omega^1 &=& 2dq_1\wedge\alpha^1-dq_2\wedge\alpha^2+(-1+p_2-q_1q_2)\alpha^1\wedge\alpha^2,\\
 \Omega^2 &=& dq_1\wedge\alpha^1-2dq_2\wedge\alpha^2+(1-p_1+q_1q_2)\alpha^1\wedge\alpha^2,\\
  \Omega^3 &=& dr_1\wedge\alpha^1+dp_2\wedge\alpha^2-(2r_1q_2+3q_1p_2)\alpha^1\wedge\alpha^2 ,\\ 
   \Omega^4 &=& dp_1\wedge\alpha^1+dr_2\wedge\alpha^2-(3p_1q_2+2r_2q_1)\alpha^1\wedge\alpha^2,\\ 
    \Omega^5 &=& dr_2\wedge \alpha^1-dr_1\wedge \alpha^2 +4(q_1r_1-q_2r_2)\alpha^1\wedge \alpha^2,\\
     \Omega^6 &=& 2w_1dp_2\wedge\alpha^2-2w_2dp_1\wedge\alpha^1+3w_3(p_2-p_1)\alpha^1\wedge\alpha^2.
      \end{array}
\right.
       \end{equation}

\begin{remark}
Notice that the differential system $(\mathcal{I},\alpha^1\wedge\alpha^2)$ is
quasi-linear.
\end{remark}

 

To discuss the involutivity of the system we compute the polar spaces
of 1-dimensio\-nal inte\-gral ele\-ments.
On $P$, we consider the coframe 
\begin{equation}
   (\alpha^1,\alpha^2,\eta^1,\dots,\eta^{16},dq_1,dq_2,dp_1,dp_2,dr_1,dr_2)
     \end{equation}
and its dual frame field
\begin{equation}
 \left(\frac{\partial}{\partial \alpha^1},\frac{\partial}{\partial
  \alpha^2},\frac{\partial}{\partial
   \eta^1},\dots,\frac{\partial}{\partial
    \eta^{16}},\frac{\partial}{\partial q_1},\frac{\partial}{\partial
     q_2},\frac{\partial}{\partial p_1},\frac{\partial}{\partial p_2},\frac{\partial}{\partial
      r_1}, \frac{\partial}{\partial r_1}\right).
      \end{equation} 
The $1$-dimensional integral elements $E$ of the system are of the form 
\begin{equation}
  E =[V(a,b,c,d)],\quad V(a,b,c,d)= a_{j}\frac{\partial}{\partial
   \alpha^j}+b_j\frac{\partial}{\partial q_j}+c_j\frac{\partial}{\partial p_j}
     +d_j\frac{\partial}{\partial r_j}.
    \end{equation}
Thus, the manifold of 1-dimensional integral elements $\mathcal{V}_1\cong P\times \RP^7$. 
A $1$-dimensional integral element is admissible if and only if $(a_1)^2+(a_2)^2\neq 0$. 
The polar equations of a given $E\in\mathcal{V}_1$ are
\begin{equation} 
 \eta^\alpha=0, \quad \alpha =1,\dots,16
  \end{equation}
and
\begin{equation}
i_{V}\Omega^\beta =0, \quad \beta =1,\dots,6,
   \end{equation} 
which read
\begin{gather*}
 2a_1dq_1-a_2dq_2 = [a_2(1-p_2+q_1q_2)+2b_1)]\alpha^1+[a_1(p_2-q_1q_2-1)-b_2]\alpha^2,\\
 a_1dq_1-2a_2dq_2 = [b_1-a_2(1-p_1+q_1q_2)]\alpha^1+[a_1(1-p_1+q_1q_2)-2b_2]\alpha^2,\\
 a_1dr_1+a_2dp_2 =[d_1 +a_2(2r_1q_2+3q_1p_2)]\alpha^1+[c_2-a_1(2r_1q_2+3q_1p_2)]\alpha^2,\\
 a_1dp_1+a_2dr_2 = [c_1 +a_2(3p_1q_2+2r_2q_1)]\alpha^1+[d_2 -a_1(3p_1q_2+2r_2q_1)]\alpha^2,\\
 a_2dr_1-a_1dr_2 =[4a_2(q_1r_1-q_2r_2) -d_2]\alpha^1+[d_1 - 4a_1(q_1r_1-q_2r_2)]\alpha^2,\\
w_2a_1dp_1 -w_1a_2dp_2 = [w_2c_1 -\frac{3}{2}a_2w_3(p_1-p_2)]\alpha^1
      +[\frac{3}{2}a_1w_3(p_1-p_2)-w_1c_2]\alpha^2.
\end{gather*}
Therefore, if $a_1a_2(w_1(a_1)^2 -w_2(a_2)^2)\neq 0$,
the polar equations are linearly independent and the dimension of the
polar space $H(E)$ of $E$ is $2$. Thus $H(E)$ is the only 2-dimensional integral element that
contains $E$.
This shows that the system is in involution and that the general
integral submanifolds depend on six functions in one variable.




\bibliographystyle{amsalpha}

\end{document}